\documentclass[preprint,sort&compress,12pt]{elsarticle}
\usepackage{comment}

\usepackage{xcolor}
 
\usepackage[algo2e]{algorithm2e}
\usepackage{amssymb}
\usepackage{amsthm}
\usepackage{amsmath}
\usepackage{mathrsfs}
\usepackage{mathtools}
\usepackage{textcomp}
\usepackage{gensymb}
\usepackage{graphicx}
\usepackage{subfig}
\usepackage{float}
\usepackage{caption}

\usepackage{tabu}
\usepackage{diagbox}
\usepackage{tikz}
\usepackage[utf8]{inputenc}
\usepackage{pgfplots}
\pgfplotsset{compat=newest}
\usepgfplotslibrary{groupplots}
\usepgfplotslibrary{dateplot}
\usepackage{tikzscale}
\usetikzlibrary{arrows,automata}
\usepackage{algorithm}
\usepackage{algpseudocode}
\usepackage{multirow}
\usepackage{enumerate}
\usepackage{array}

\usepackage[flushleft]{threeparttable}
\usepackage{lineno}
\usepackage{fullpage}
\usepackage{xcolor}
\usepackage{bbm}
\usepackage{listings}
\usepackage[colorlinks=true]{hyperref}
\usepackage{url}
\usepackage{siunitx}
\graphicspath{ {./figs/} }
\biboptions{comma,square}
\graphicspath{ {./figs/} }
\journal{Elsevier}
\newcommand{\norm}[1]{\left\lVert#1\right\rVert}

\theoremstyle{definition}

\begin{document}
\begin{frontmatter}
\hypersetup{pdfauthor=author}

 \title{A Bi-fidelity Ensemble Kalman Method for PDE-Constrained Inverse Problems in Computational Mechanics}



\author[ndAME,ndCICS]{Han Gao}
\author[ndAME,ndCICS]{Jian-Xun Wang\corref{corxh}}

\address[ndAME]{Department of Aerospace and Mechanical Engineering, University of Notre Dame, Notre Dame, IN}
\address[ndCICS]{Center for Informatics and Computational Science, University of Notre Dame, Notre Dame, IN}
\cortext[corxh]{Corresponding author. Tel: +1 574-631-5302}
\ead{jwang33@nd.edu}

\begin{abstract}
Mathematical modeling and simulation of complex physical systems based on partial differential equations (PDEs) have been widely used in engineering and industrial applications. To enable reliable predictions, it is crucial yet challenging to calibrate the model by inferring unknown parameters/fields (e.g., boundary conditions, mechanical properties, and operating parameters) from sparse and noisy measurements, which is known as a PDE-constrained inverse problem. In this work, we develop a novel bi-fidelity (BF) ensemble Kalman inversion method to tackle this challenge, leveraging the accuracy of high-fidelity models and the efficiency of low-fidelity models. The core concept is to build a BF model with a limited number of high-fidelity samples for efficient forward propagations in the iterative ensemble Kalman inversion. Compared to existing inversion techniques, salient features of the proposed methods can be summarized as follow: (1) achieving the accuracy of high-fidelity models but at the cost of low-fidelity models, (2) being robust and derivative-free, and (3) being code non-intrusive, enabling ease of deployment for different applications. The proposed method has been assessed by three inverse problems that are relevant to fluid dynamics, including both parameter estimation and field inversion. The numerical results demonstrate the excellent performance of the proposed BF ensemble Kalman inversion approach, which drastically outperforms the standard Kalman inversion in terms of efficiency and accuracy. 
\end{abstract}

\end{frontmatter}


\section{Introduction}
\label{sec:intro}
%

Partial differential equation (PDE)-based mathematical models have been pervasively used in modeling complex physical systems (e.g., complex fluid flows \cite{catalano2003numerical,iaccarino2003immersed,kalitzin2005near,evans2013isogeometric1,evans2013isogeometric2,luo2008immersed}) for various engineering applications \cite{evans2009n,evans2013isogeometric3,zheng2009computational}. With ever-increased computational resource and advances in the theoretical understanding of the physics, a model tends to be increasingly sophisticated, which, however, poses great challenges on calibrating considerable unknown/uncertain model parameters (e.g., boundary conditions, mechanical properties, or other control parameters). Model calibration (i.e., parameter estimation) is an \emph{inverse problem}, which can be formulated as,
\begin{equation}
\label{eq:inverse}
\mathrm{Find}\ \mathbf{z} \in \mathcal{I}_{z} \ \mathrm{given \ observations\ } \mathbf{y} = \mathcal{M}(\mathbf{z}) + \boldsymbol{\sigma_d} \in \mathcal{I}_y, 
\end{equation}
where $\mathcal{M}$ is the parameter-to-observable mapping that maps the parameter vector $\mathbf{z}\in \mathcal{I}_{z}$ to the observation space $\mathcal{I}_y$, and such a mapping usually involves a PDE-based physical model, e.g., computational fluid dynamics (CFD) models for fluid dynamics; $\boldsymbol{\sigma_d}$ represents process errors that reflect measurement uncertainty and model inadequacy. 

Traditionally, the inverse problem defined in Eq.~\ref{eq:inverse} is solved by minimizing the mismatch between model prediction $\mathcal{M}(\mathbf{z})$ and observation $\mathbf{y}$ in a least-squares manner, which is known as the \emph{variational approach}. To efficiently solve the variational optimization problem, local gradients of the cost function with respect to the parameters are needed and such derivative information is often obtained via adjoint method~\cite{plessix2006review,giles2000introduction}. Although the adjoint-based variational method has been successfully applied for a variety of PDE-constrained inverse problems~\cite{dow2011quantification,talnikar2017unsteady,pires2001tsunami,wang2019discrete,rangarajan2020adjoint,wang2020identification}, it has two major limitations. \underline{First},  the development of an adjoint solver is highly code-intrusive, which requires substantial effort to implement, especially for existing large-scale legacy codes of computational mechanics~\cite{nielsen2010discrete}. \underline{Second}, variational optimization is often solved deterministically and thus the prediction uncertainty associated with observation noises and model-form errors cannot be quantified~\cite{asch2016data}. As an alternative,  the inverse problem (Eq.~\ref{eq:inverse}) can also be posted in a probabilistic way, which is to search for the posterior distribution $p(\mathbf{z}|\mathbf{y})$ instead of the optimal values of inferred quantities $\mathbf{z}$ given observation data $\mathbf{y}$ based on Bayes' theorem, referred to as the \emph{Bayesian approach}~\cite{cotter2009bayesian,ou2019new,zhang2019efficient}. The advantage of the Bayesian formulation is that the inversion tends to be more robust and associated uncertainties can be naturally quantified through statistical properties \cite{zhang2018quantification,zhang2018effect}. In general, the posterior is obtained by Monte Carlo (MC)-based sampling techniques such as importance sampling and Markov chain Monte Carlo (MCMC) methods~\cite{li2015adaptive,morzfeld2018iterative,uzun2019structural}, where a large number of samples are required and thus the computation becomes intractable for any nontrivial forward model (e.g., CFD simulations of complex flows). 

Iterative ensemble Kalman methods (IEnKM) have been recently developed for efficiently solving inverse problems in a Bayesian manner~\cite{iglesias2013ensemble,liu2016nonlocal,sun2015statistical,liu2016determining,sun2015hybrid,sun2013identification,feng2015simultaneous,sun2016probabilistic}. IEnKM can be viewed as the application of ensemble Kalman filtering (EnKF), which has been demonstrated to be highly successful for dynamic state estimation in data assimilation applications~\cite{evensen2003ensemble,evensen2009data,carrassi2018data}, towards steady inverse problems~\cite{schillings2017analysis}. The basic idea of IEnKM is to iteratively perform Kalman analysis for parameter estimation in a nonlinear or non-Gaussian setting, which was originally proposed in refs~\cite{chen2012ensemble,chen2013levenberg,iglesias2013ensemble} and has been recently analyzed~\cite{schillings2017analysis,schillings2018convergence,evensen2018analysis} and improved in terms of convergence~\cite{blomker2019well,wu2019improving} and identifiability~\cite{iglesias2016regularizing,wu2019adding,zhang2019regularization,albers2019ensemble,chada2019incorporation}. Compared to the adjoint-based variational methods or sampling-based Bayesian approaches, the salient advantages of IEnKM are that (i) it is \emph{derivative-free} and \emph{code non-intrusive}; (ii) it can provide reliable estimations with a \emph{small size ensemble} (e.g., $\mathcal{O}(10)$ of ensemble members). Many successful applications of IEnKM have been witnessed in recent years for various inverse problems of computational mechanics, including field inversion in complex fluid flows~\cite{iglesias2015iterative,wang2016data,wang2018inferring,tang2018tsuflind}, calibration of turbulence models~\cite{kato2013approach,xiao2016quantifying,wu2016bayesian,mons2016reconstruction}, and state-parameter estimation in physiological systems~\cite{arnold2017uncertainty,lal2017non,wang2019data}, and others~\cite{zoccarato2016data,iglesias2018ensemble,sousa2019computational,yang2020non}. 

In contrast to sampling-based Bayesian approaches like MCMC, the ensemble in IEnKM is only utilized to estimate the state/parameter covariance, which is in lieu of the gradient to enable a derivative-free optimizer in a variational manner~\cite{schillings2017analysis}. It has been demonstrated that a small number of ensemble members are sufficient for data assimilation (i.e., state estimation) in many noisily observed dynamical systems~\cite{evensen2009data}. Nonetheless, when it comes to steady-state inverse problems with strong nonlinearity, multiple iterations of the Kalman analysis are required to further refine the mean estimation~\cite{iglesias2013ensemble}. As a result, a considerable number of forward model evaluations are still needed (i.e., the number of iterations $\times$ ensemble size), introducing a substantial computational burden, particularly when the forward simulation is expensive, e.g., three-dimensional (3-D) CFD simulations for complex flows. In order to overcome this obstacle in dealing with large-scale inverse problems, developing a more efficient forward propagation scheme is of great significance. One promising strategy is to leverage lower fidelity forward models that are less expensive and less accurate (e.g., coarse-mesh and/or un-converged solutions). By optimally combining models with varying levels of accuracy and cost, the computational cost of multiple forward evaluations is expected to be largely reduced, which is known as a multi-fidelity method. In this paper, we will develop a novel bi-fidelity iterative ensemble Kalman method (BF-IEnKM) for efficient parameter/field inversion, by leveraging the accuracy of high-fidelity models and efficiency of low-fidelity models. Specifically, a bi-fidelity surrogate model is constructed for the forward propagation based on the recently proposed multi-fidelity stochastic collocation scheme~\cite{narayan2014stochastic,zhu2014computational}. The parameter ensemble is iteratively propagated to the observables using the constructed bi-fidelity surrogate and updated by the Kalman analysis. The performance of the proposed method is evaluated on a number of inverse problems in fluid dynamics applications. The most relevant work to this paper is the multilevel ensemble Kalman filter (MLEnKF) proposed by Hoel et al.~\cite{hoel2016multilevel}, where a multilevel Monte Carlo (MLMC) sampling strategy is embedded into the Monte Carlo step of the standard EnKF. Contrary to MLMC, the high-fidelity solutions in this work will be utilized as basis functions for a low-rank approximation of the solution space, while the low-fidelity model will be used to inform global searches over the parameter space and assist the high-fidelity solution reconstruction. To the knowledge of the authors there has yet to be an extension of the bi-fidelity method to the ensemble-based Kalman inversion.

The rest of the paper is organized as follows. In Section~\ref{sec:meth}, the methodology and algorithm of the proposed BF-IEnKM are introduced. To demonstrate the effectiveness, several PDE-constrained inverse problems (including parameter estimation and field inversion) are studied in Section~\ref{sec:result}, where the proposed BF-IEnKM is compared with the standard IEnKM in terms of efficiency and accuracy. Finally, Section~\ref{sec:conclusion} concludes the paper.

\section{Methodology \label{sec:meth}}

\subsection{Overview}
\begin{figure}[htp]
    \centering
    \includegraphics[width=0.9\textwidth]{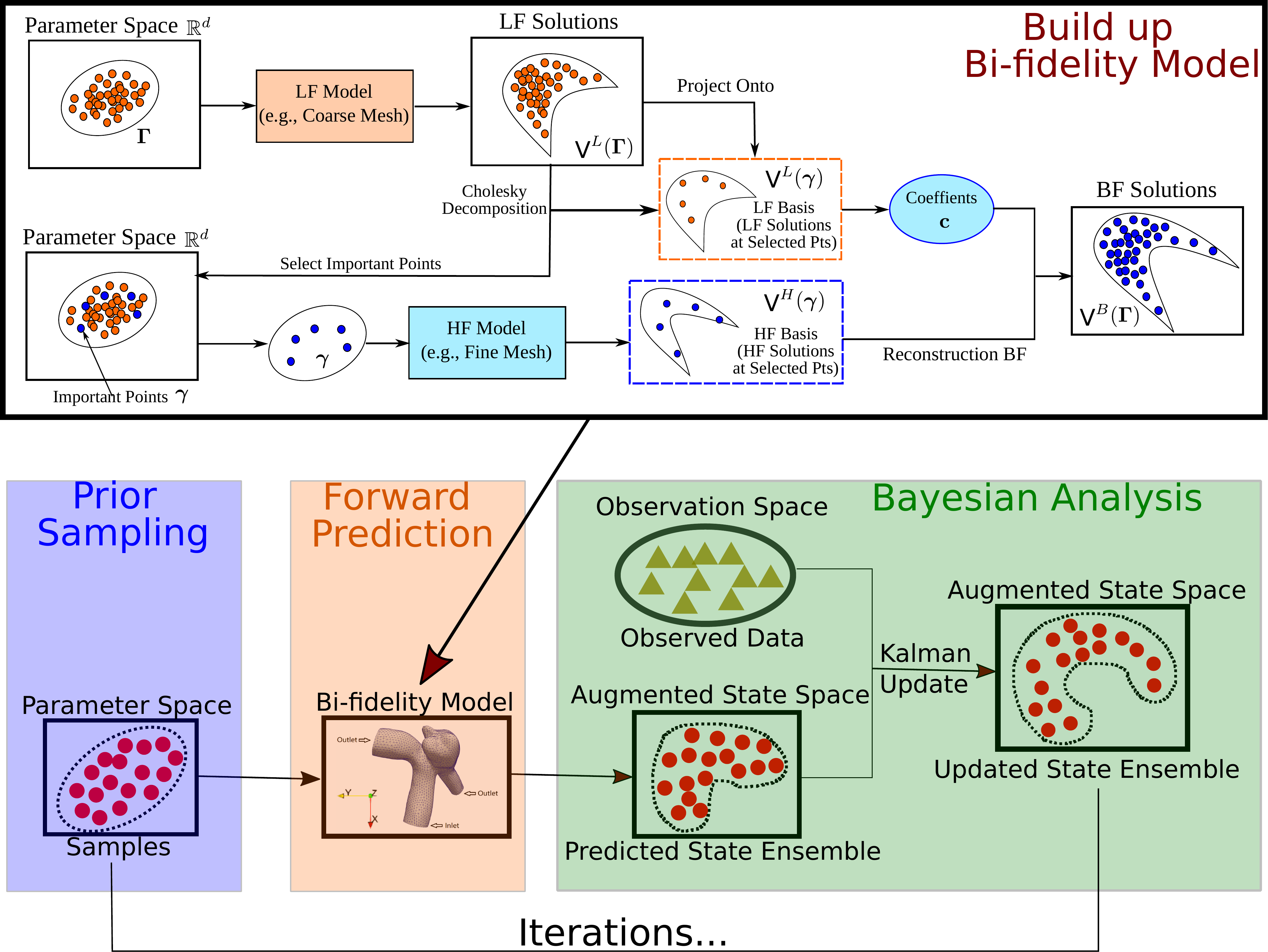}
    \caption{Schematic of the proposed bi-fidelity iterative ensemble Kalman methods (BF-IEnKM) for PDE-constrained inverse problems.}
    \label{fig:EnKFMethod}
\end{figure}

The goal here is to solve the inverse problem defined by Eq.~\ref{eq:inverse}, where the parameter-to-observable mapping $\mathcal{M}$ is usually composed of a forward solution operator $\mathcal{F}: \mathcal{I}_{z} \to \mathcal{I}_{x}$ that maps the parameter space $\mathcal{I}_{z}$ to the state space $\mathcal{I}_{x}$ and an observation operator $\mathcal{H}: \mathcal{I}_{x} \to \mathcal{I}_y$ that projects the simulated state $\mathbf{x} \in \mathcal{I}_{x}$ to the observation space $\mathcal{I}_y$. If measurements $\mathbf{y} \in \mathcal{I}_y$ are given, the unknown/uncertain parameters $\mathbf{z}$ can be inferred by iteratively applying Bayesian analysis (Kalman update) as shown in Fig.~\ref{fig:EnKFMethod}.
Specifically, (1) an ensemble of parameters (e.g., boundary conditions) are sampled from the prior distribution; (2) all these parameter samples are propagated through the forward mapping $\mathcal{F}$ to obtain an ensemble of predicted states (e.g., velocity/pressure fields); (3) the state ensemble is augmented with parameters, and the augmented state-parameter ensemble is then updated by assimilating observation data $\mathbf{y}$ based on the Kalman analysis formula, which will be detailed later. The steps (2) and (3) are iteratively conducted until a certain convergence criterion is satisfied. As mentioned above, in these procedures, the forward propagation step can be computationally burdensome since the functional mapping $\mathcal{F}$ is usually defined by a system of PDEs, which has to be solved numerically. Namely, the PDEs are spatially (and temporally) discretized into a large-scale finite-dimensional algebraic system, and the solution accuracy primarily relies on sufficient mesh resolution and numerical iterations, which is thus time-consuming. On the other hand, coarsening the mesh grid or reducing the number of numerical iterations can largely reduce the computational cost but significantly sacrifice the solution accuracy. Here we refer to the fine-mesh converged solution as high-fidelity (HF) model while the coarse-mesh unconverged solution as low-fidelity (LH) model. To enable an efficient and reliable forward propagation, a bi-fidelity (BF) strategy leveraging advantages of both LF and HF models will be applied to develop a bi-fidelity iterative ensemble Kalman method (BF-IEnKM).  Figure~\ref{fig:EnKFMethod} shows the overall schematics of the proposed BF-IEnKM, and the details of each part, such as the BF propagation scheme and iterative ensemble Kalman inversion algorithm will be presented in following subsections. 


\subsection{Bi-fidelity forward propagation}
The essence of the proposed BF-IEnKM is to use the bi-fidelity stochastic collocation strategy~\cite{narayan2012stochastic,narayan2014stochastic,zhu2014computational,gao2019bi,hampton2018practical,skinner2019reduced,fairbanks2020bi} to forward propagate the ensemble from the parameter space ($\mathcal{I}_{z} \subseteq \mathbb{R}^d$) to the solution space ($\mathcal{I}_{x} \subseteq \mathbb{R}^n$). The key is to construct a low-rank approximation of the solution space based on a limited number of HF solutions $\mathbf{v}^H(\gamma)$ on a small set of ``important" parameter points $\gamma \subset \mathcal{I}_{z}$. The selection of these ``important" points is informed by exploring the parameter space using a large number of cheap low-fidelity simulations. The BF forward propagation can be formulated either in an online or offline manner. In the online formulation, a new BF model is constructed for every Kalman iteration, while the offline BF surrogate is pre-trained and remains unchanged for all Kalman iterations. Although the former can adapt the changes due to Kalman updates, the training cost is significantly higher than that of the latter one. Considering the moderate ensemble size in IEnKM, we choose to construct the BF forward solver in an offline manner. Namely, the BF propagation scheme has two phases: (1) \underline{offline training phase}, where a low-rank approximation of HF solution space is pre-trained based on both HF and LH simulation data, and (2) \underline{online prediction phase}, where only the LF simulation is needed to propagate a new point in the parameter space. The number of HF training data will be estimated based on an empirical error estimate proposed in \cite{gao2019bi}. 


\subsubsection{Offline training phase}
To select the small set of important parameter points, a large number ($M\geqslant 1$) of cheap LF simulations are conducted on a prescribed nodal set $\Gamma = \{\mathbf{z}_1,...,\mathbf{z}_M\}$ covering the entire parameter space $\mathcal{I}_z$. A subset $\gamma^{(m)} = \{\mathbf{z}_1,...,\mathbf{z}_m\} \subset \Gamma$ of ``important" points are selected from $\Gamma$ using a greedy algorithm proposed in~\cite{narayan2012stochastic,narayan2014stochastic}, where a new point is iteratively selected such that the LF solution $\mathbf{v}^L(\mathbf{z}_k)$ of the newly added point $\mathbf{z}_k$ has the furthest distance from the space spanned by the LF solutions of the existing selected parameters $\gamma^{(k-1)} = \{\mathbf{z}_1,...,\mathbf{z}_{k-1}\}$. This can be formulated as an optimization problem,
\begin{subequations}
	\label{eqn:maxmizer}
	\begin{alignat}{2}
	\mathbf{z}_{k} = \arg\max_{\mathbf{z}\in\Gamma} d\bigg(\mathbf{v}^L(\mathbf{z}), \mathbb{U}^L(\gamma^{(k-1)})\bigg), \\
	\gamma^{(k)} = \gamma^{(k-1)}\cup\{\mathbf{z}_k\},
	\end{alignat}
\end{subequations}
where $d(\boldsymbol{v},W)$ is the distance function between the vector $\boldsymbol{v}\in \mathbf{v}^L(\Gamma)$ and subspace $W\subset \mathbb{U}^L(\Gamma)=\textrm{span}\{\mathbf{v}^L(\mathbf{z}_1),\cdots,\mathbf{v}^L(\mathbf{z}_M)\}$. This optimization can be solved by factorizing the Gramian matrix $G$ of LF solutions  $V^L(\Gamma)=[\mathbf{v}^L(\mathbf{z}_1),..,\mathbf{v}^L(\mathbf{z}_M)]^T$ based on the pivoted Cholesky decomposition~\cite{zhu2014computational},
\begin{equation}
G = P^TLL^TP,
\end{equation}
where $L$ is a low-triangular matrix; $P$ is a permutation matrix such that $P\Gamma$ provides an order of ``importance", and the index corresponding to the first $m$ columns can be used to identify the important parameter points for HF simulations. The implementation is based on the Algorithm 1 in ref~\cite{gao2019bi} (See \ref{Appendix:ImportantPointSelection}). Once the small set of important points $\gamma^{(m)}$ are selected, HF solutions $V^H(\gamma^{(m)})=[\mathbf{v}^H(\mathbf{z}_1),..,\mathbf{v}^H(\mathbf{z}_m)]^T$ on $\gamma^{(m)}$ can be utilized as the basis functions for a low-rank approximation of the solution space. 

\subsubsection{Online prediction phase}
The BF solution of a new parameter point $\mathbf{z} \in \mathcal{I}_z$ is constructed based on the HF solution basis functions $V^H(\gamma^{(m)})$ obtained in the offline training phase. To compute the reconstruction coefficients, a LF-based approximation is used by assuming that HF and LF reconstructions share the same coefficients. Namely, only the LF solution $\mathbf{v}^L(\mathbf{z})$ of the new point $\mathbf{z}$ is simulated and projected onto the LF solution basis functions $V^L(\gamma^{(m)})$ to find the projection coefficients $\mathbf{c}^L(\mathbf{z})$. In summary, the BF propagation can be performed as follows:
\begin{enumerate}
	\item Perform LF simulation at $\mathbf{z}$,
	\begin{equation}	 
	\mathbf{v}^L(\mathbf{z})=\mathcal{F}^{L}(\mathbf{z}).
	\end{equation}	
	\item Project the LF solution $\mathbf{v}^L$ onto the pre-computed LF basis $V^L(\gamma^{(m)})$ to compute reconstruction coefficients,
	\begin{equation}	 
	\mathbf{c}^L=((G^L)^{-1})^T\mathbf{v}^L(\mathbf{z}).
	\end{equation}
	\item Reconstruct the BF solution using reconstruction coefficients $\mathbf{c}^L$ and pre-computed HF basis $V^L(\gamma^{(m)})$ as, 
	\begin{equation}
	\label{eq:reconstruct}
	\mathbf{v}^H(\mathbf{z}) \approx \mathbf{v}^B(\mathbf{z})=\sum_{k=1}^m c_k\mathbf{v}^H(\mathbf{z}_k)
	\end{equation}
\end{enumerate} 

\subsubsection{A priori error estimation}
In order to evaluate if the LF model is sufficiently informative and determine how many HF simulations are necessary for the desired accuracy, \emph{a priori} assessment of the LF model quality and error estimation of the BF predictions are practically useful~\cite{narayan2014stochastic,hampton2018practical}. To this end, an empirical error estimation approach proposed in~\cite{gao2019bi} is employed here. Firstly, two non-dimensional scalar metrics $R_s(\mathbf{z})$ and $R_e(\mathbf{z})$ are introduced. The model similarity metric $R_s(\mathbf{z})$ measures the similarity between LF and HF models, defined as,
\begin{equation}
\label{eqn:relativeDist}
R_{s}(\mathbf{z}) = \frac{d(\mathbf{v}^H(\mathbf{z}),\mathbb{U}^H(\gamma^{k}))}{||\mathbf{v}^H(\mathbf{z})||}/\frac{d(\mathbf{v}^L(\mathbf{z}),\mathbb{U}^L(\gamma^{k}))}{||\mathbf{v}^L(\mathbf{z})||}.
\end{equation}
An informative (i.e., good) LF model leads to $R_{s} \approx 1$. The error component ratio $R_e(\mathbf{z})$ defined as, 
\begin{equation}
\label{eqn:rError}
R_{e}(\mathbf{z}) = \frac{|| P_{\mathbb{U}^H({\gamma^k})}\mathbf{v}^H(\mathbf{z}) -  \mathbf{v}^B(\mathbf{z})||}{d(\mathbf{z},\mathbb{U}^H(\gamma^{k}))},
\end{equation}
which measures the balance between the in-plane error $|| P_{\mathbb{U}^H({\gamma^k})}\mathbf{v}^H(\mathbf{z}) -  \mathbf{v}^B(\mathbf{z})||$ and the relative distance $d(\mathbf{z},\mathbb{U}^H(\gamma^{k}))$. When $R_e$ becomes a large value (e.g., $10^1$), indicating that the in-plane error is dominant over the distance error, we should stop generating new HF samples because the greedy search algorithm is based on the distance component. The error estimate of the BF solution at any new point $\mathbf{z}_*$ is given by,
\begin{equation}
\label{eq:eb3}
\frac{||\mathbf{v}^H(\mathbf{z}_*) - \mathbf{v}^B(\mathbf{z}_*)||}{||\mathbf{v}^H(\mathbf{z}_*)||}
\leq c \max_{\mathbf{z}\in\Gamma}(\frac{d(\mathbf{v}^L(\mathbf{z}),\mathbb{U}^L(\gamma^{k}))}{||\mathbf{v}^L(\mathbf{z})||}) (1+R_{e}(\mathbf{z}_{k+1})),
\end{equation}
where constant $c$ is set to be 1 as recommended in \cite{gao2019bi}. The error estimate only requires the LF data and pre-selected HF data and it has been demonstrated effective when the LF model is informative ($R_s \approx 1$) and the in-plane error is not dominant ($R_e < 10^1$)~\cite{gao2019bi}.

\subsection{Bi-fidelity iterative ensemble Kalman inversion}
Similar to the original iterative ensemble Kalman method (IEnKM), the initial parameter ensemble is sampled from the prior distribution and then propagated to the solution (state) ensemble via the BF forward propagation scheme described above. Both the predicted state and parameter samples are corrected based on Kalman analysis, and the updated parameter ensemble becomes the new prior in the next iteration. 

\subsubsection{Prior sampling step}
Initially, an parameter ensemble $\gamma^{(n_s)}_0 = \{\mathbf{z}_i\}_{i=1}^{n_s}$ is drawn from the prior distribution $p(\mathbf{z})$ using the Latin hypercube sampling (LHS) method \cite{stein1987large, shields2016generalization}, where $n_s$ is the size of the ensemble. 

\subsubsection{Forward prediction step}
At iteration $k$, the parameter ensemble  $\gamma^{(n_s)}_k$ is first propagated via the LF model,
\begin{equation}
V^{L}_{k}(\gamma^{(n_s)}_k) = \mathcal{F}^{L}(\gamma^{(n_s)}_k).
\end{equation}
Then the LF solution ensemble $V^{L}_{k}$ is projected onto the LF basis $V^L(\gamma^{(m)})$ pre-computed in the offline phase to obtain the reconstruction coefficients,
\begin{equation}
\mathbf{c}_k= ((G^L)^{-1})^T V^{L}_{k}.
\end{equation}
The coefficients $\mathbf{c}_k=[c_{i}^1,...,c_{i}^m]$ coupled with the pre-computed HF basis $V^H(\gamma^{(m)})$ are used to reconstruct the BF solution ensemble $V^{B}_{k}(\gamma^{(n_s)}_k)$ based on Eq.~\ref{eq:reconstruct}.



\subsubsection{Bayesian analysis step}
In the Bayesian analysis step, observation data $\mathbf{y}$ will be incorporated to update the parameters based on their correlations. To obtain the correlations between the parameter $\mathbf{z}$ and the state $\mathbf{v}$, an augmented state vector is defined as $\mathbf{x}=[\mathbf{v},\mathbf{z}]$. The augmented state ensemble propagated via the BF forward model at iteration $k$ is denoted by $X_k^{B} = [V_k^B, \gamma_k^{(n_s)}]$. The indirect and noisy observation data $\mathbf{y}$ can be expressed as,
\begin{equation}
\mathbf{y}=\mathcal{H}(\tilde{\mathbf{v}})+\boldsymbol{\sigma}_d
\end{equation}  
where $\tilde{\mathbf{v}}$ represent the true state and $\boldsymbol{\sigma}_d$ is modeled as an zero-mean Gaussian noise with covariance matrix $P_d$. The nonlinear projection operator $\mathcal{H}$ can be easily converted to a linear projection matrix $H$ by augmenting the observables into the state. The augmented state is then corrected using the Kalman analysis formula,
\begin{equation}
\hat{X}_k^B=X_k^B+P_m^k H^T (H P_m^k H^T + P_d)^{-1}(\mathbf{y}-H X_k^B),
\end{equation}
where $P_m^k$ is the covariance matrix of the augmented state ensemble $X_k^B$. Repeat the prediction and analysis steps until the ensemble is statistically converged (e.g., the data mismatch is much lower than the data noise level).

\section{Numerical Results}
\label{sec:result}
In order to illustrate the accuracy, efficiency, and applicability of the proposed bi-fidelity iterative ensemble Kalman method (BF-IEnKM), we present numerical experiments of three different inverse problems constrained by PDEs. The first example (Case 1) considers parameter estimation in a two-dimensional (2-D) convection-diffusion system with weak discontinuity. The second example (Case 2) is a model calibration problem in Reynolds-averaged Navier-Stokes (RANS) turbulence modeling, where the coefficients of a RANS turbulence model is reversely determined from sparse and noisy mean flow measurements. In the third example (Case 3), a more challenging inverse problem is studied, where the inflow velocity field of a 3-D patient-specific aneurysm flow is inferred given sparse and noisy internal flow data. In all three cases, synthetic observation data are used. Namely, sparsely sampled simulation results with specified parameter/boundary conditions (i.e., ground truth) are used as data, where artificial Gaussian noises can be added. The HF forward propagation is performed on fine mesh grids, while the LF model is constructed by significantly downsampling the mesh. The number of Kalman iterations is determined based on the stopping criterion in~\cite{wang2019data}. To assess the accuracy of inferred parameters $\mathbf{z}^{p}$, the relative error metric ($e$) is defined as,
\begin{equation}
e=\sqrt{\frac{\norm{\mathbf{z}^{p}-\tilde{\mathbf{z}}}_{L^2}}{\norm{\tilde{\mathbf{z}}}_{L^2}}},
\end{equation}
where $\tilde{\mathbf{z}}$ is the ground truth, $\mathbf{z}^{p}$ represents the posterior mean of the parameter ensemble, and $\norm{\cdot}_{L^2}$ is the L2 norm defined in the parameter space. The forward solvers (with both HF and LF discretization) are implemented on an open-source finite volume method (FVM) platform, OpenFOAM, where spatial derivatives were discretized with the second-order central scheme for both convection and diffusion terms. A second-order implicit time-integration scheme was used to discretize the temporal derivatives. For the incompressible fluid problems considered in cases 2 and 3, the continuity and moment equations were solved using the SIMPLE (semi-implicit method for pressure-linked equations) algorithm~\cite{caretto1973two}, and collocated grids with Rhie and Chow interpolation were used to prevent the pressure-velocity decoupling~\cite{rhie1983numerical}. The case settings are summarized in Table~\ref{tab:ParamterSummary}.
\begin{table}[htp]
	\centering
	\begin{tabular}{|c|c|c|c|}
		\hline
		\multicolumn{1}{|c|}{\diagbox{Control parameters}{Cases}}&
		\multicolumn{1}{|c|}{Case 1}&
		\multicolumn{1}{|c|}{Case 2}&
		\multicolumn{1}{|c|}{Case 3}\\
		\hline
		\multicolumn{1}{|c|}{Mesh grids (HF)}&
		\multicolumn{1}{|c|}{10,000}&
		\multicolumn{1}{|c|}{12,225}&
		\multicolumn{1}{|c|}{51,648}\\
		\hline
		\multicolumn{1}{|c|}{Mesh grids (LF)}&
		\multicolumn{1}{|c|}{49}&
		\multicolumn{1}{|c|}{75}&
		\multicolumn{1}{|c|}{4464}\\
		\hline
		\multicolumn{1}{|c|}{$\#$ of (LF/HF) simulations for BF}&
		\multicolumn{1}{|c|}{$1000/15$}&
		\multicolumn{1}{|c|}{$1000/20$}&
		\multicolumn{1}{|c|}{$2000/60$}\\
		\hline
		\multicolumn{1}{|c|}{Parameters to be inferred ($\mathbf{z}$)}&
		\multicolumn{1}{|c|}{$D_T$}&
		\multicolumn{1}{|c|}{$C_{1\varepsilon}, C_{2\varepsilon}$}&
		\multicolumn{1}{|c|}{$\mathbf{u}(\mathbf{x}), \mathbf{x}\in\partial D_{in}$}\\		
		\hline
		\multicolumn{1}{|c|}{$\#$ of data (percentage of HF grids)}&
		\multicolumn{1}{|c|}{$2\%$}&
		\multicolumn{1}{|c|}{$1\%$}&
		\multicolumn{1}{|c|}{$0.3\%$}\\
		\hline
		\multicolumn{1}{|c|}{Data noise (Std)}&
		\multicolumn{1}{|c|}{$0\%$ of truth}&
		\multicolumn{1}{|c|}{$40\%$ of truth}&
		\multicolumn{1}{|c|}{$40\%$ of truth}\\
		\hline
		\multicolumn{1}{|c|}{Minimum Kalman iterations}&
		\multicolumn{1}{|c|}{3}&
		\multicolumn{1}{|c|}{2}&
		\multicolumn{1}{|c|}{4}\\
		\hline
		\multicolumn{1}{|c|}{Ensemble size}&
		\multicolumn{1}{|c|}{30}&
		\multicolumn{1}{|c|}{25}&
		\multicolumn{1}{|c|}{200}\\
		\hline
	\end{tabular}
	\caption{Computational parameters in BF-IEnKM for three test cases}
	\label{tab:ParamterSummary}
\end{table}

\subsection{Case 1: parameter estimation in 2-D convection-diffusion system}
In the first example, a 2-D convection-diffusion equation is considered,
\begin{equation}
\frac{\partial T}{\partial t}+ \nabla\cdot(\mathbf{u} T )- \nabla^2(D_T T)=0,
\end{equation}
where $T$ is the scalar concentration field, $\mathbf{u}$ is the velocity field, and $D_T$ is the diffusion coefficient. For simplicity, both $\mathbf{u}$ and $D_T$ are set to be constants over the entire domain. The convection-diffusion equation describes mass transport of a scalar quantities $T$ (e.g., solute species) based on the combined effect of both convection and diffusion. When the diffusivity $D_T$ is small, the convection effect is dominant and sharp gradients (or even shock) will present in the steady-state $T$ field, which is notoriously difficult to resolve by classic numerical methods, e.g., finite volume or continuous Galerkin methods \cite{canuto2019uncertainty}. In general, a very fine mesh can be used to capture the large gradients, which, however, will significantly increase the computational cost and pose great challenges on inverse problems. To better evaluate the proposed BF-IEnKM, the goal here is to reversely determine the diffusion coefficient $D_T$, whose true value is small, leading to a nearly discontinuous concentration field. 

\begin{figure}[htp]
	\centering
	\subfloat[Prior mean]{\includegraphics[width=0.4\textwidth]{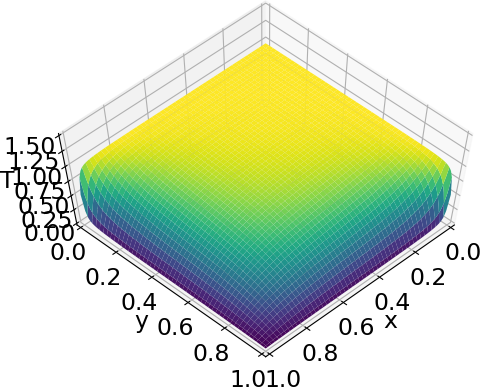}}
	\subfloat[Posterior mean (LF-IEnKM)]{\includegraphics[width=0.4\textwidth]{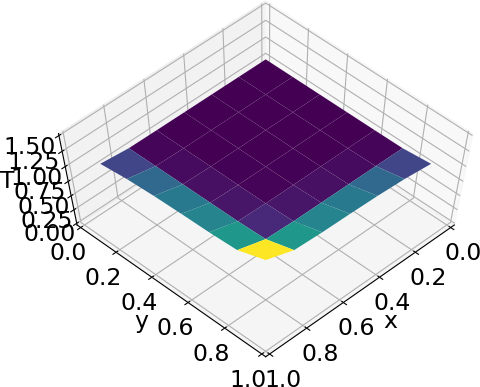}}
	\vfill
	\subfloat[Posterior mean (BF-IEnKM)]{\includegraphics[width=0.4\textwidth]{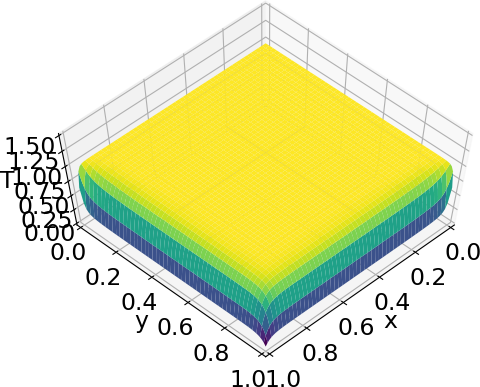}}
	\subfloat[Ground truth]{\includegraphics[width=0.4\textwidth]{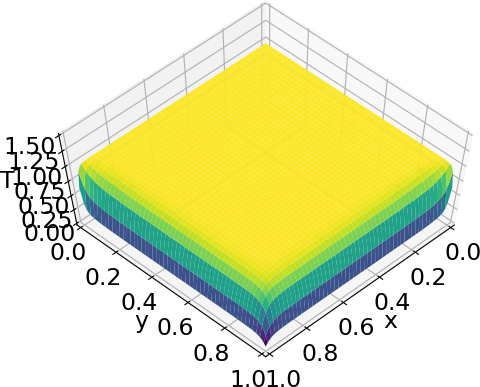}}
	\caption{Scalar concentration $T$ field with (a) prior $D_T$, (b) posterior $D_T$ obtained by LF-IEnKM, (c) posterior $D_T$ obtained by BF-IEnKM, and (d) true $D_T$.}
	\label{fig:case1TContour}
\end{figure}

Consider a 2-D domain $[0,1]\times[0,1]$ with a constant convection velocity $\mathbf{u} = [1, 1]^T$. Dirichlet boundary condition (BC) are specified, i.e., $T = 1$ on the left and bottom sides ($[0,1]\times\{0\}$, $\{0\}\times[0,1]$) and $T = 0$ on the right and top sides ($[0,1]\times\{1\}$, $\{1\}\times[0,1]$). The true diffusion parameter is set as $\tilde{D}_T = 0.025$, leading to sharp gradients in the steady-state $T$ field, as shown in Fig.~\ref{fig:case1TContour}d. A HF model with $100\times100$ grids and a LF model with only $7\times7$ grids are used to build the BF forward surrogate. In the offline phase, the solutions of $1,000$ parameters randomly drawn from a uniform distribution of $D_T \in [0.02, 1.98]$ are explored by the LH model and $15$ of them are selected as important parameter points where the HF simulations are conducted for the BF reconstruction. The synthetic observation data are obtained by randomly sampling $2\%$ of the HF solution field $\tilde{T}$ with the true $\tilde{D}_T$ (no noise is added in this case). The initial $D_T$ ensemble of 30 members is sampled from a uniform prior $U(0.15,0.25)$ and are updated by the BF-IEnKM. 
\begin{figure}[htp]
	\centering
	\includegraphics[width=0.6\textwidth,height=0.28\textwidth]{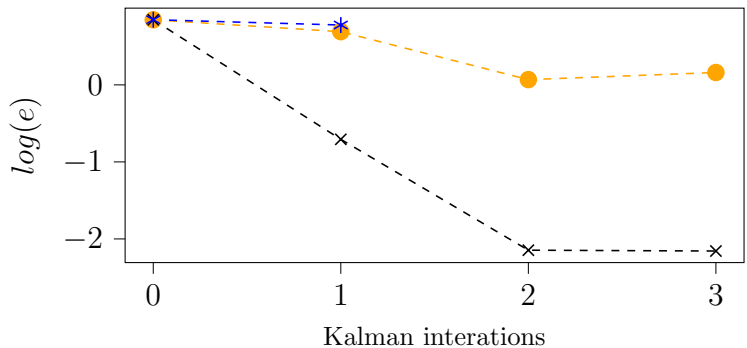}
	\caption{Comparison of inferred diffusion coefficient $D_T$ using LF, BF, and HF models with noise-free observations. BF sample mean (\ref{plot:case1ParaBF}), LF sample mean (\ref{plot:case1ParaLF}), and the HF sample mean (\ref{plot:case1ParaHF}). } 
	\label{fig:case1Para0}
\end{figure}
As shown in Fig.~\ref{fig:case1Para0}, only with three Kalman iterations, the relative error of BF-inferred $D_T$ (posterior sample mean $D_T^{BF}=0.02483$) has been reduced less than $1\%$ (\ref{plot:case1ParaBF}), while the error of the LF-inferred result ($D_T^{LF}=0.0611$) remains large, which is over $100\%$ (\ref{plot:case1ParaLF}). Note that after the offline BF reconstruction, the online cost of the BF-based inversion is almost the same as that of the LF-based inversion. It's also interesting to compare the performance of the BF-IEnKM with that of the purely HF-based inversion using the \emph{same amount of computational budget}. Namely, we only can afford one Kalman iteration with 15 HF samples, which costs the same as the BF-based inversion does, including the budget for both offline BF training and online inversion. However, the purely HF-based inversion simply fails and the relative error of the posterior mean ($D_T^{HF}=-0.1250$) was barely reduced, which is above $100\%$ (\ref{plot:case1ParaHF}). This is because both the iterations and ensemble size are not sufficiently large.

The posterior (inferred) $D_T$ ensemble can be propagated to the concentration $T$ field to further evaluate the peformance. Figure~\ref{fig:case1TContour} compared the contours of posterior mean fields of $T$ by LF-based inversion (Fig.~\ref{fig:case1TContour}b) and BF-based inversion (Fig.~\ref{fig:case1TContour}c), benchmarked against the ground truth (Fig.~\ref{fig:case1TContour}d), where the prior mean (Fig.~\ref{fig:case1TContour}a) is also plotted for comparison.
It can be seen that the LF posterior mean becomes even worse compared to the prior and fails to give a physical result because of numerical issues induced by insufficient mesh resolution. In contrast, the BF posterior mean agrees with the ground truth very well, and the sharp gradients near the two boundary sides can be accurately captured, demonstrating the effectiveness of the proposed approach.

\subsection{Case 2: RANS turbulence model calibration for 2-D separated flow}
\begin{figure}[htp]
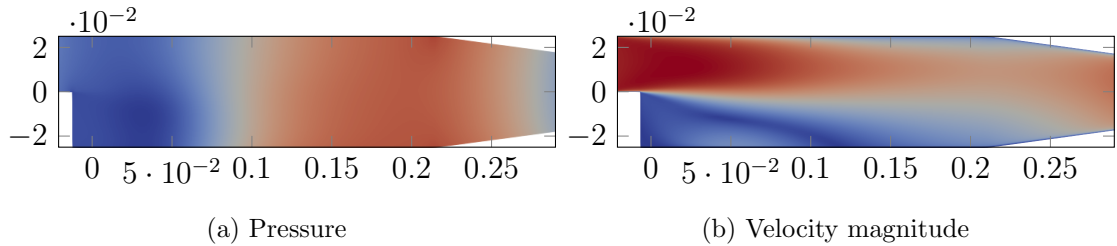

	\centering
	\subfloat[Pressure]
	{\includegraphics[width=0.45\textwidth,height=0.15\textwidth]{./pGT.tikz}}
	\subfloat[Velocity magnitude]
	{\includegraphics[width=0.45\textwidth,height=0.15\textwidth]{./UGT.tikz}}
	\caption{Mean pressure and velocity magnitude fields of the back ward facing step at $Re=25000$ with the $k$-$varepsilon$ model}.
	\label{fig:case2Contour}
\end{figure}

\begin{figure}[htp]
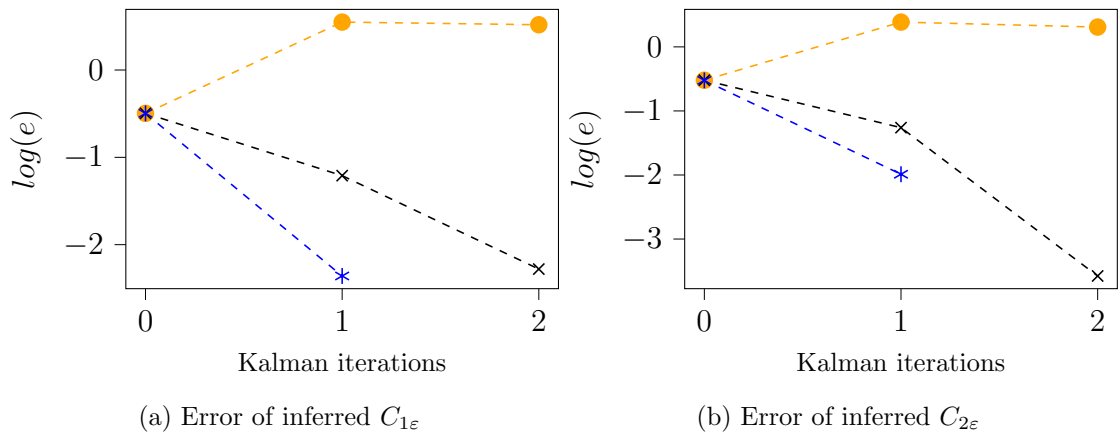

    \centering
    \subfloat[Error of inferred $C_{1\varepsilon}$]
    {\includegraphics[width=0.45\textwidth,height=0.3\textwidth]{./parameter_case2_0.tikz}}
    \subfloat[Error of inferred $C_{2\varepsilon}$]
    {\includegraphics[width=0.45\textwidth,height=0.3\textwidth]{./parameter_case2_1.tikz}}
    \caption{Comparison of inferred $C_{1\varepsilon}$ and $C_{2\varepsilon}$ using LF, BF, and HF models with noisy observations ($40\%$ Gaussian noise). BF sample mean (\ref{plot:case1ParaBF}), LF sample mean (\ref{plot:case1ParaLF}), and the HF sample mean (\ref{plot:case1ParaHF}).}
    \label{fig:case2Para}
\end{figure}

Turbulent flow is ubiquitous in natural and industrial processes, which can be simulated using CFD techniques. Despite the increasing computational power, high-fidelity turbulence simulations, e.g., direct numerical simulation or large eddy simulation, are still not affordable for most engineering problems. Numerical models based on Reynolds-Averaged Navier-Stokes (RANS) equations are still the workhorse tools for practical applications. To close the RANS equations, many turbulence closures (e.g., $k$-$\varepsilon$, $k$-$\omega$, and Spalart-Allmaras models) have been proposed in the past decades based on certain specific assumptions~\cite{durbin2018some}, which, however, are not universally applicable and have large model-form errors particularly for flows with curvature, strong pressure gradient, and massive separations~\cite{spalart2009detached}. There are growing interests in quantifying and reducing the model-form uncertainty in RANS turbulence closures based on Bayesian inference~\cite{wu2016bayesian,xiao2016quantifying,wang2016incorporating,kato2013approach,edeling2018bayesian} and physics-informed machine learning~\cite{wang2017physics,ling2016machine,duraisamy2019turbulence,wang2019prediction,yang2019predictive}. Among them, calibrating the parameters of a given closure model using sparse flow measurements is one strategy to improve the RANS modeling~\cite{kato2013approach,edeling2018bayesian}. In this example, we test our proposed BF inversion approach to calibrate a RANS model. To be specific, two coefficients $C_{1\varepsilon}$ and $C_{2\varepsilon}$ of the standard $k$-$\varepsilon$ model~\cite{launder1983numerical} will be reversely determined based on sparse and noisy mean flow measurements. The turbulent flow behind a backward-facing step with a Reynolds number of $Re = 25,000$~\cite{fureby2010homogenization} is studied. The ground truth of the two coefficients are set as $\tilde{C}_{1\varepsilon} = 1.44$ and $\tilde{C}_{2\varepsilon} = 1.92$, which are the recommended values in~\cite{launder1983numerical}. The corresponding mean velocity and pressure fields with $\tilde{C}_{1\varepsilon}$ and $\tilde{C}_{2\varepsilon} $ are shown in Figs.~\ref{fig:case2Contour}a and~\ref{fig:case2Contour}b, respectively.
The HF model has 12,225 mesh grids, while the LF one only has 75 cells. The synthetic data are obtained by observing the true mean flow on $1\%$ randomly selected HF grids, and corrupted with $40\%$ Gaussian noises. To reconstruct the BF surrogate, 1,000 LF samples are generated by uniformly sampling the parameter space $(C_{1\varepsilon},C_{2\epsilon})\in[1,2]\times[1.8, 3.0]$, and HF simulations are conducted on the first 20 most important parameters as the basis. For the online Kalman iteration, an ensemble of 25 members is drawn from the prior and updated by the synthetic observations. The iterations histories of the inferred $C_{1\varepsilon}$ and $C_{2\varepsilon}$ are plotted in Figs.~\ref{fig:case2Para}a and~\ref{fig:case2Para}b, respectively. 
The purely LF-based inversion simply fails and the inferred results ($C_{1\varepsilon}^{LF}=6.1902$ and $C_{2\varepsilon}^{LF}=5.8339$) become even worse compared to the initial guess (the relative error is higher than $100\%$ \ref{plot:case2ParaLF}). The performance of the BF-IEnKM are excellent as the inferred values ($C_{1\varepsilon}^{BF}=1.4476$ and $C_{2\varepsilon}^{BF}=1.9205$) are very close to the ground truth ($\tilde{C}_{1\varepsilon} = 1.44$ and $\tilde{C}_{2\varepsilon} = 1.92$) in only two iterations (relative errors are lower than $1\%$ \ref{plot:case2ParaBF}). Similarly, the purely HF-based inversion is conducted using the same computational budget of the BF-based inversion, which affords one Kalman iteration with 20 HF samples. It can be seen that the purely HF-based inversion yields satisfactory results (\ref{plot:case2ParaHF}) compared to the LF results. However, the accuracy is lower than that of the BI-IEnKM method. The HF-inferred $C_{1\varepsilon}$ is nearly as the same accurate as the BF-inferred one, while for $C_{2\varepsilon}$, the BF-inferred result is over one magnitude more accurate than the HF-inferred result ($C_{1\varepsilon}^{HF}=1.4337$ and $C_{2\varepsilon}^{HF}=1.9003$).

\subsection{Case 3: field inversion in patient-specific cerebral aneurysm flow}
\begin{figure}[htp]
	\centering
	{\includegraphics[height=0.08\textwidth]{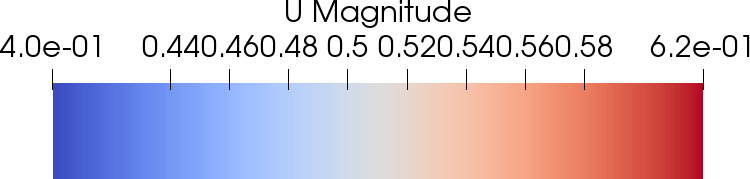}}
	\subfloat[Prior]
	{\includegraphics[width=0.25\textwidth]{./IGInlet.tikz}}
	\subfloat[LF posterior]
	{\includegraphics[width=0.25\textwidth]{./LFInlet.tikz}}
	\subfloat[BF posterior]
	{\includegraphics[width=0.25\textwidth]{./BFInlet.tikz}}
	\subfloat[Ground Truth]
	{\includegraphics[width=0.25\textwidth]{./GTInlet.tikz}}
	\caption{Contour comparison of the (a) prior mean, (b) LF posterior mean, and (c) BF posterior mean of the inlet velocity field, benchmarked against (d) the ground truth.}
	\label{fig:case3_inlet40noise}
\end{figure}

In the third case, the proposed BF-IEnKM is applied to solve a field inversion problem in a cardiovascular flow application. Consider blood flows in a 3-D patient-specific cerebral aneurysm bifurcation~\cite{valen2018real}, which can be simulated based on the steady incompressible Navier-Stokes equations under assumptions of rigid walls and Newtonian fluids. The reliability of model predictions (e.g., velocity, pressure, and wall shear stress) largely depends on the boundary conditions, e.g., inflow velocity field, which are usually uncertain or unknown. Therefore, the goal here is to infer a complex 3D inlet velocity field (including non-zero secondary flow) from sparse and noisy internal flow measurements. To specify the prior distribution of velocity inlet, a stationary zero-mean Gaussian random field $f(\mathbf{x})$ is introduced to model the prior uncertainty in the inlet:
\begin{equation}
f(\mathbf{x})\sim \mathcal{GP}(0,K(\mathbf{x},\mathbf{x}')),\;\; 
K(\mathbf{x},\mathbf{x}')=\sigma_0^2\exp{(\frac{|\mathbf{x}-\mathbf{x}'|}{2l^2})},
\end{equation}
where $K(\mathbf{x},\mathbf{x}')$ is the exponential kernel function, $\sigma_0$ and $l$ defined the standard deviation and length scale of the random field, respectively. The random field is expressed in a compact form using Karhunen-Loeve (K-L) expansion \cite{tipping1999probabilistic},
\begin{equation}
f(\mathbf{x})=\sum_{i=1}^{n_k\to \infty}\sqrt{\lambda_i}\phi_i(\mathbf{x})\omega_i,
\end{equation}
where $\lambda_i$ and $\phi_i(\mathbf{x})$ are eigenvalues and eigenfunctions of the kernel $K$; $\omega_i$ are independent and identically distributed  (i.i.d.) random variables with zero mean and unit variance. We further truncate the KL expansion with a finite number ($n_k=3$) of basis to approximate the stochastic process. In this case, a Gaussian random field with $l=2\times10^{-3}$m and $\sigma_0=0.1$m/s is added to each component of the inlet velocity with the first three K-L modes, capturing over $90\%$ energy. The prior distribution of the inlet velocity fields can be written as
\begin{equation}
\mathbf{u}_{in}=\mathbf{u}_{in}^{base}+[f_x(\mathbf{x}),f_y(\mathbf{x}),f_z(\mathbf{x})],
\end{equation}
where the base streamwise velocity magnitude $\|\mathbf{u}_{in}^{base}\|$is set as $0.509$ m/s ($Re = 345$). It is noted that the inlet velocity has variations not only in streamwise direction but also in the in-plane directions, which induces complex secondary flows. Inferring the inflow velocity field can be formulated as a nine-dimensional inverse problem with the following parameters to be inferred,
\begin{equation}
\mathbf{z} = [ 
\underbrace{\omega_1, \omega_2, \omega_3}_{\mathrm{x\ streamwise}},
\underbrace{\omega_4, \omega_5, \omega_6}_{\mathrm{y\ in-plane}},
\underbrace{\omega_7, \omega_8, \omega_9}_{\mathrm{z\ in-plane}}
]^T \in \mathcal{I}_z \subseteq \mathbb{R}^9.
\label{eqn:paraspacecase3}
\end{equation}
One realization is set as the ground truth, and the corresponding inflow is shown in Fig~\ref{fig:case3_inlet40noise}d.

\begin{figure}[htp]
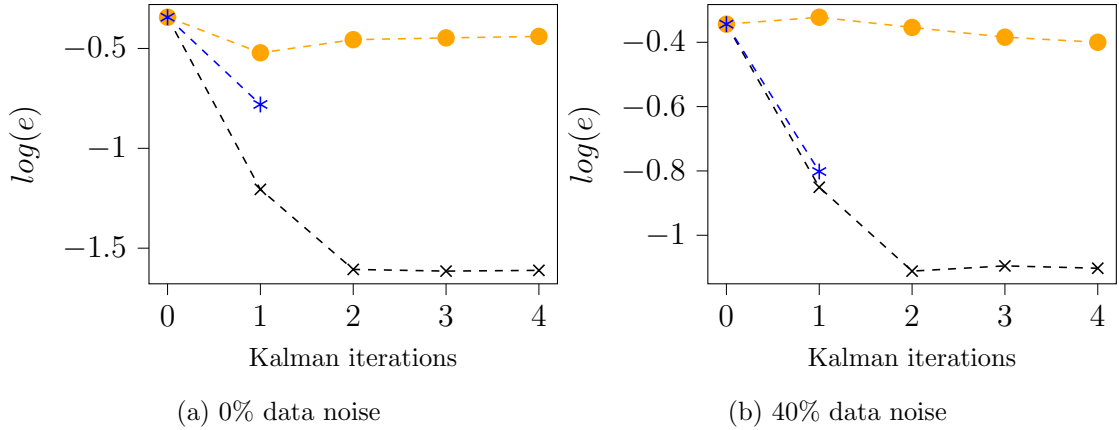

	\centering
	\subfloat[$0\%$ data noise]
	{{\includegraphics[width=0.45\textwidth,height=0.3\textwidth]{./inletErr_0Noise.tikz}}}
	\subfloat[$40\%$ data noise]
	{{\includegraphics[width=0.45\textwidth,height=0.3\textwidth]{./inletErr.tikz}}}
	\caption{Relative error of inferred inlet velocity fields ($\mathbf{u}_{in}$) by LF, BF, and HF models with (a) noise-free velocity data, and (b) velocity data with $40\%$ Gaussian noises. BF sample mean (\ref{plot:case3ErrorBF}), LF sample mean (\ref{plot:case3ErrorLF}), and the HF mean (\ref{plot:case3ErrorHF}).}
	\label{fig:case3error}
\end{figure}

The BF model is built upon an HF model with 51,648 mesh grids and an LF model with 4,464 mesh grids. In the offline phase, the total number of 2,000 LF samples are drawn from a nine-dimensional multivariate Gaussian distribution, where 60 of them are selected for the HF simulations. In the online Kalman inversion, an initial ensemble of 200 parameters are drawn from the prior and iteratively updated by assimilating synthetic observations, which are obtained by randomly sampling $0.3\%$ of true internal flow velocity results. Figure~\ref{fig:case3error} shows the relative error of the inferred inflow velocity field v.s. iterations.  
When noise-free data is assimilated (Fig.~\ref{fig:case3error}a), the relative error of the BF-inferred $\mathbf{u}_{in}$ (\ref{plot:case3ErrorBF}) can be reduced to about $2.5\%$ after four iterations, while the LF-inferred result (\ref{plot:case3ErrorLF}) is barely improved compared to the prior. Even when the synthetic data are perturbed with a $40\%$ Gaussian noise (Fig.~\ref{fig:case3error}a), the performance of the proposed BF-IEnKF is still impressive and the accuracy of the BF-inferred $\mathbf{u}_{in}$ is above $92\%$. With the same computational budget, only 60 HF samples with one iteration can be afforded to perform the purely HF-based inversion, and the performance (\ref{plot:case3ErrorHF}) is worse than the BF-based approach for both cases with noise-free or noisy data. To clearly demonstrate the effectiveness of the proposed methods, the contour plots of inferred velocity inlets (Case with $40\%$ noisy data) are visualized in Fig.~\ref{fig:case3_inlet40noise}.
The mean field of the prior inflow ensemble (Fig.~\ref{fig:case3_inlet40noise}a) has a uniform streamwise velocity component and no secondary flows since zero-mean Gaussian processes are specified. The true inlet field (Fig.~\ref{fig:case3_inlet40noise}d) presents a complex flow pattern, where the velocity magnitude is non-uniformly distributed and in-plane secondary flows can be founded. It is clear that the BF-inferred inflow velocity field (Fig.~\ref{fig:case3_inlet40noise}c) agrees with the ground truth very well, which cannot be captured by the purely LF-based approach (Fig.~\ref{fig:case3_inlet40noise}b). Note that the successful inference is achieved by assimilating highly sparse and noisy data ($0.3\%$ of velocity information corrupted with $40\%$ noise), demonstrating the merit and effectiveness of the proposed BF-based iterative ensemble Kalman method.

\begin{figure}[htp]
	\centering
	\subfloat[Prior (flow)]
	{\includegraphics[width=0.25\textwidth]{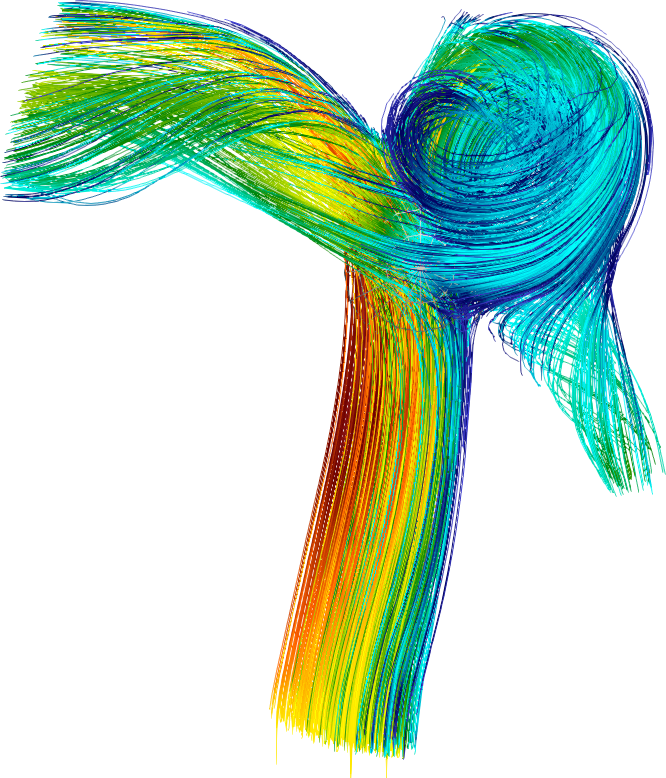}}
	\hfill
	\subfloat[BF posterior (flow)]
	{\includegraphics[width=0.25\textwidth]{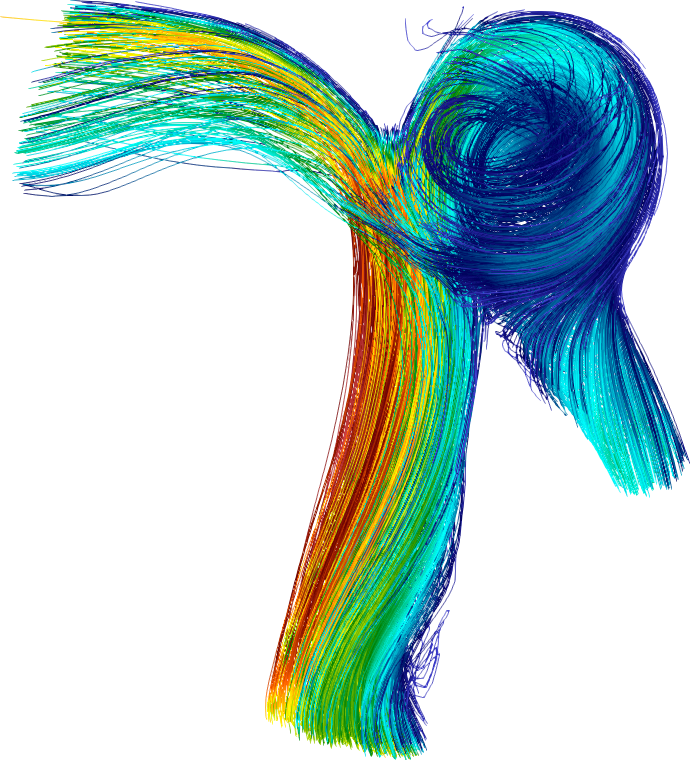}}
	\hfill
	\subfloat[Truth (flow)]
	{\includegraphics[width=0.25\textwidth]{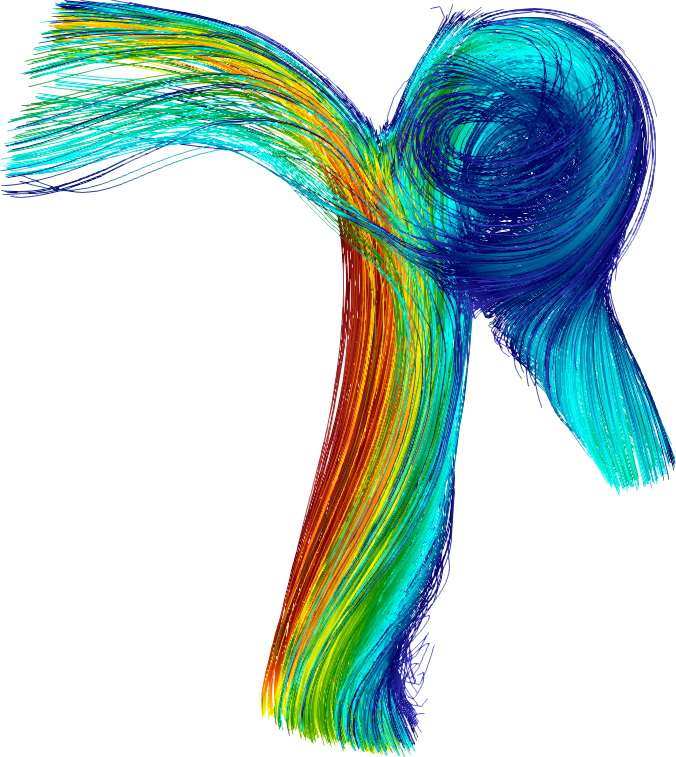}}\\	
	\subfloat[Prior (pressure)]
	{\includegraphics[width=0.25\textwidth]{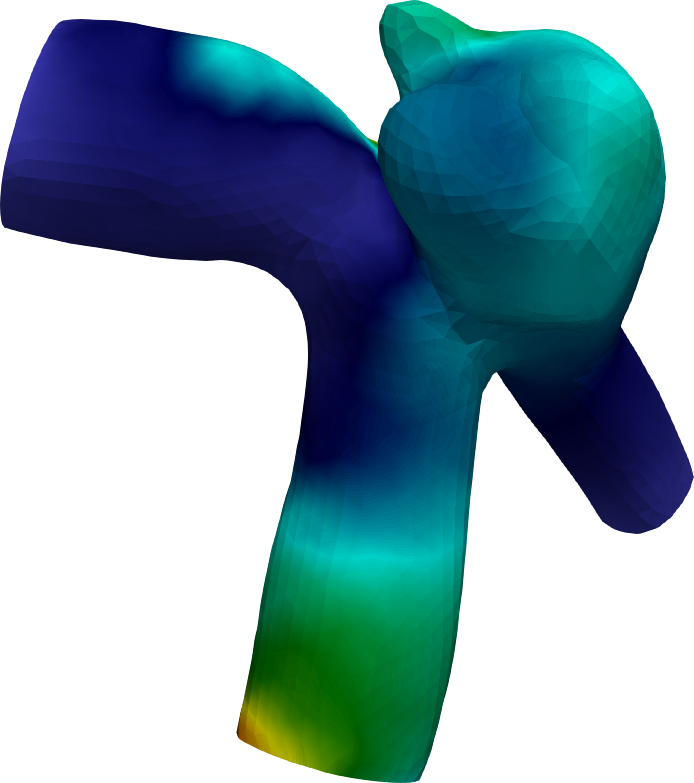}}
	\hfill
	\subfloat[BF posterior (pressure)]
	{\includegraphics[width=0.25\textwidth]{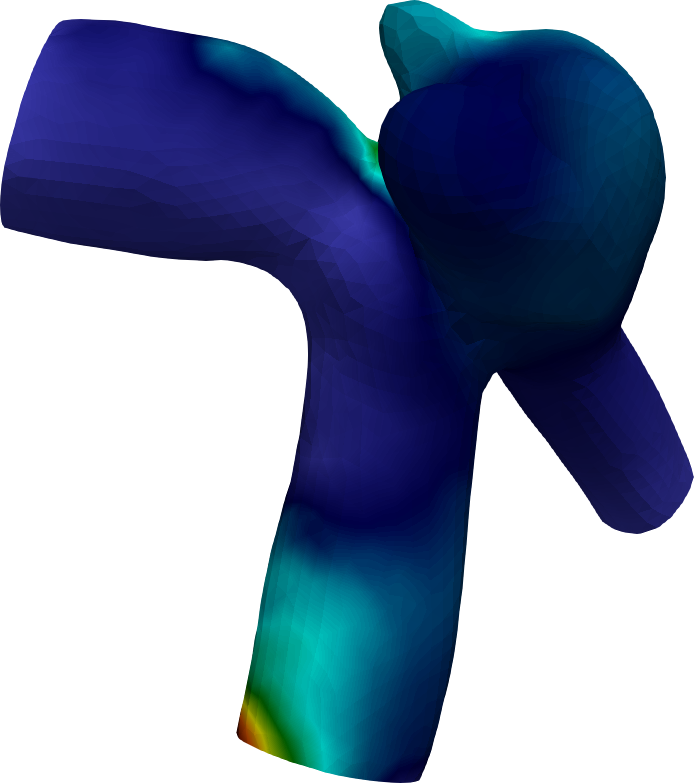}}
	\hfill
	\subfloat[Truth (pressure)]
	{\includegraphics[width=0.25\textwidth]{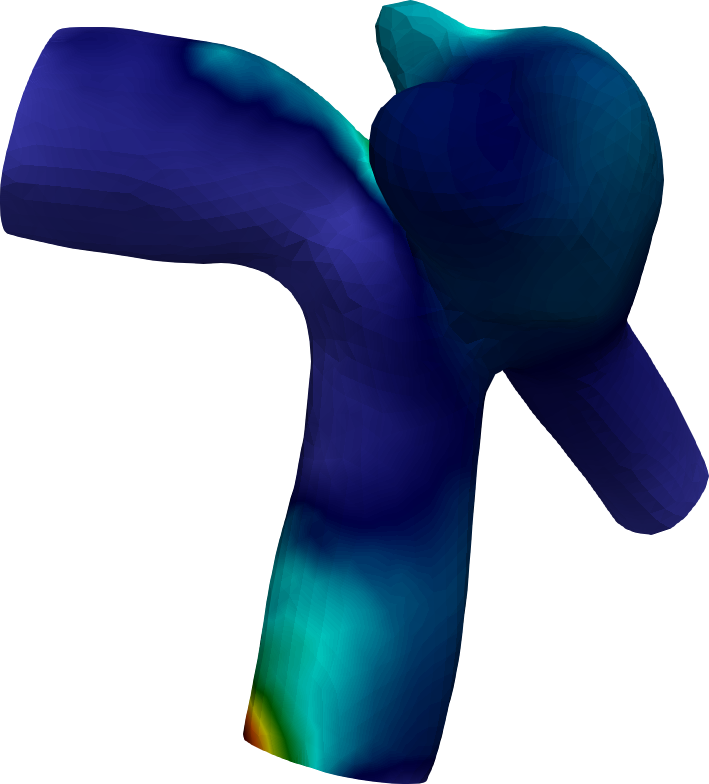}}\\
	\subfloat[Prior (WSS)]
	{\includegraphics[width=0.25\textwidth]{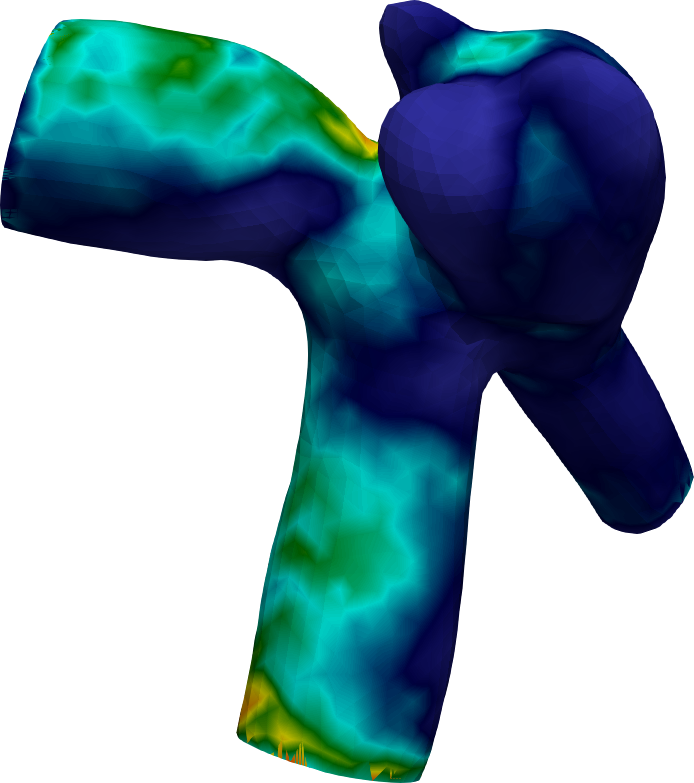}}
	\hfill
	\subfloat[BF posterior (WSS)]
	{\includegraphics[width=0.25\textwidth]{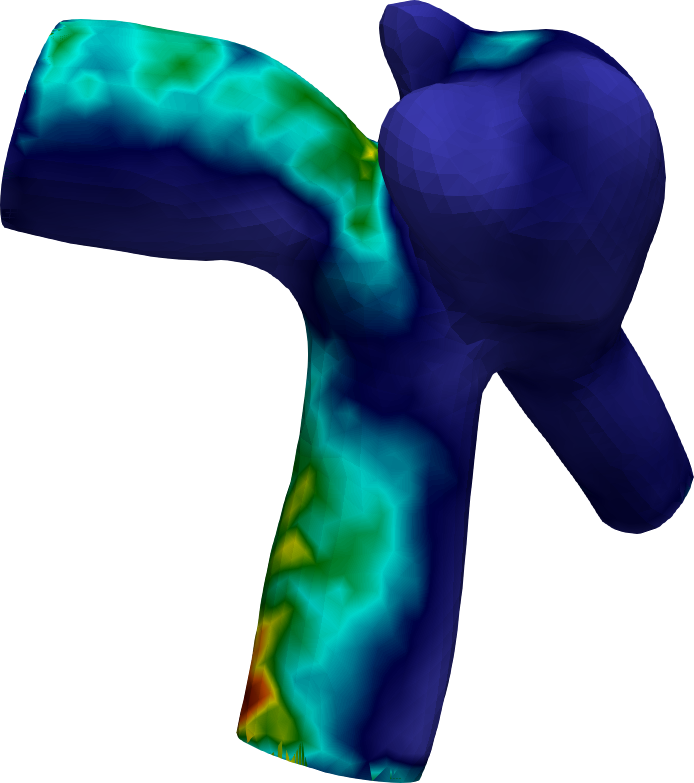}}
	\hfill
	\subfloat[Truth (WSS)]
	{\includegraphics[width=0.25\textwidth]{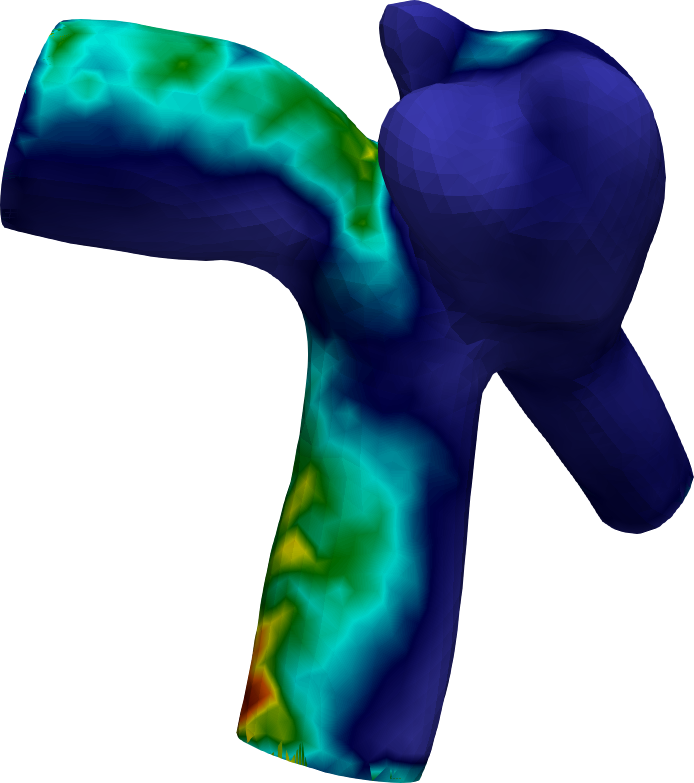}}	
	\caption{Propagated full-field flow information, including (a, b, c) flow streamline, (d, e, f) surface pressure, and (g, h, i) wall shear stress distribution.}
	\label{fig:case3_pressure}
\end{figure}
The inferred inlet velocity field can be propagated by solving the incompressible Navier-Stokes equations to obtain the full-field information of the internal flows. To further reduce the computational cost, the BF surrogate is used to enable a fast forward propagation. Namely, the posterior ensemble of the inlets are propagated via the LF solver, and the high-resolution internal velocity field, surface pressure field, and wall shear stress (WSS) field are recovered based on the BF reconstruction, which are shown in Fig.\ref{fig:case3_pressure}. Because of the secondary flows of the inlet, strong streamline distortions can be seen in the main vessel, aneurysm bulb, and two bifurcation arms, which are underestimated in the prior (Fig.\ref{fig:case3_pressure}a). Although only the noisy velocity information observed on $0.3\%$ of the grids are used for inference, the full-field internal velocity over the entire domain can be well recovered by the BF-based approach (Fig.\ref{fig:case3_pressure}b), which is almost identical to the ground truth (Fig.\ref{fig:case3_pressure}c). For the pressure distribution over the vessel surface, the BF prediction (Fig.\ref{fig:case3_pressure}e) also well agrees with the truth (Fig.\ref{fig:case3_pressure}f), where the high-pressure regions at left corner of the main vessel caused by the flow impingement are accurately captured. Lastly, we evaluate the WSS prediction, which is one of the most critical quantities of interest in aneurysm flows. Although the overall pattern of the prior mean WSS (Fig.\ref{fig:case3_pressure}g) is similar to the ground truth (Fig.\ref{fig:case3_pressure}i), the locations of extreme values, which are known to be highly relevant to aneurysm rupture~\cite{jou2008wall,gao2019bi}, are mispredicted. Nonetheless, our proposed BF-based approach can accurately predicted these regions, where the maximum and minimum WSS are located (Fig.\ref{fig:case3_pressure}h).

\subsection{Efficiency and Accuracy of BF-IEnKM}
Lastly, the comparison between our proposed BF-IEnKM and the standard IEnKM based on either LF or HF models are summarized in terms of the efficiency and accuracy.

\subsubsection{Speedup in online inversion phase}
Once the BF model is constructed in offline training phase, the speedup of our proposed BF-IEnKM for the online inversion directly scales up with the ensemble size, the number of iterations, and the ratio the HF/LF model evaluation costs. Table~\ref{tab:caseAllsingle} summarizes the computational cost (CPU wall time) of a single model evaluation for the BF (LF) and HF models in all three cases studied above. It can be seen the speedup is remarkable even only for one model evaluation in all cases. For instance, the single-sample speedup in the 2-D convection-diffusion problem (Case 1) is more than $4000$ times. Note that the gain on efficiency would be huge if a large ensemble is used and more iterations are conducted.
\begin{table}[H]
	\begin{center}
		\begin{tabular}{ |c|c|c|c| } 
			\hline
			\multicolumn{1}{|c|}{\diagbox{Forward model}{Case}}&\multicolumn{1}{|c|}{2-D conv-diff} &  \multicolumn{1}{|c|}{2-D turbulent flow} & \multicolumn{1}{|c|}{3-D vascular flow}\\
			\hline
			\multicolumn{1}{|c|}{BF (or LF) model}&\multicolumn{1}{|c|}{0.0130 sec} &  \multicolumn{1}{|c|}{0.097 sec} & \multicolumn{1}{|c|}{2.160 sec}\\
			\hline
			\multicolumn{1}{|c|}{HF model}&\multicolumn{1}{|c|}{57.801 sec} &  \multicolumn{1}{|c|}{45.033 sec} & \multicolumn{1}{|c|}{282.220 sec}\\
			\hline
			\multicolumn{1}{|c|}{Speedup}&\multicolumn{1}{|c|}{4446.23} &  \multicolumn{1}{|c|}{464.25} & \multicolumn{1}{|c|}{130.65}\\
			\hline
		\end{tabular}
		\caption{Computational cost (CPU wall time,) of a single forward evaluation of LF, BF, and HF models for all three test cases\label{tab:caseAllsingle}}
	\end{center}
\end{table}

\subsubsection{Accuracy comparison with same computational budget}
Although the proposed BF-IEnKM has been demonstrated to be remarkably efficient in the online phase, additional computational costs are still need for the offline training. 
\begin{table}[H]
	\footnotesize
	\begin{center}
		\begin{tabular}{ |c|c|c|c|c|c| } 
			\hline
			\multicolumn{1}{|c|}{Case name}&\multicolumn{1}{|c|}{2-D conv-diff} &  \multicolumn{2}{|c|}{2-D turbulent flow} & \multicolumn{2}{|c|}{3-D laminar flow}\\
			\hline
			\multicolumn{1}{|c|}{Parameters}&\multicolumn{1}{|c|}{$D_T$ (0 noise)}&\multicolumn{1}{|c|}{$C_1$ ($40\%$ noise)}&\multicolumn{1}{|c|}{$C_2$ ($40\%$ noise)}&\multicolumn{1}{|c|}{$\mathbf{u}_{in}$ (0 noise)}&\multicolumn{1}{|c|}{$\mathbf{u}_{in}$ ($40\%$ noise)}\\
			\hline
			\multicolumn{1}{|c|}{$e$ (HF)}&\multicolumn{1}{|c|}{$10^{0.78}$}&\multicolumn{1}{|c|}{$10^{-2.35}$}&\multicolumn{1}{|c|}{$10^{-1.98}$}&\multicolumn{1}{|c|}{$10^{-0.78}$}&\multicolumn{1}{|c|}{$10^{-0.80}$}\\
			\hline
			\multicolumn{1}{|c|}{$e$ (BF)}&\multicolumn{1}{|c|}{$10^{-2.16}$}&\multicolumn{1}{|c|}{$10^{-2.28}$}&\multicolumn{1}{|c|}{$10^{-3.58}$}&\multicolumn{1}{|c|}{$10^{-1.61}$}&\multicolumn{1}{|c|}{$10^{-1.10}$}\\
			\hline
			\multicolumn{1}{|c|}{Budget}&\multicolumn{1}{|c|}{975 s} &  \multicolumn{2}{|c|}{326 s} & \multicolumn{2}{|c|}{12,600 s}\\
			\hline
			\multicolumn{1}{|c|}{Device name}&\multicolumn{5}{|c|}{Intel® Xeon(R) Gold 6138 CPU @ 2.00GHz × 40}\\
			\hline
		\end{tabular}
		\caption{Accuracy comparison between BF-IEnKM and HF-IEnKM with the same computational budget for all three test cases\label{tab:caseAllBudget}}
	\end{center}
\end{table}\vspace{-2em}
To more fairly compare the proposed method with the standard single-fidelity (i.e., HF-based) IEnKM, the computation costs for both offline training and online inversion of the BF-IEnKM will be considered. Namely, the total computational budget will be fixed, and the accuracy of the two approaches are compared. The comparison results are summarized in Table~\ref{tab:caseAllBudget}. For all three cases, with same computational budget, the BF-IEnKM predictions are much more accurate than those of the HF-based IEnKM, regardless of the data noise level. For most inferred parameters, the prediction accuracy of the BF-based inversion is one order of magnitude higher than that of the standard Kalman inversion, showing the merits of our proposed method.

\section{Conclusion}
\label{sec:conclusion}
This paper presents a novel bi-fidelity iterative ensemble Kalman method (BF-IEnKM) for solving PDE-constrained inverse problem, where the advantages of high fidelity (HF) and low fidelity (LF) models are leveraged. Specifically, a large number of LF simulations are conducted to explore the parameter space and Chelosky decomposition of the LF snapshot matrix is performed to determine the important parameter nodes, on which the HF simulations are conducted. These HF solutions are used as the basis functions for massive forward propagations in the ensemble-based Kalman inversion. Therefore, the computational cost of the BF-based ensemble Kalman iterations can be significantly reduced but the predictive accuracy remains high. Moreover, the proposed method inherits the non-intrusive nature of the original iterative Kalman inversion method and thus can be easily applied to any existing numerical solvers for various applications. Several numerical experiments were conducted for solving inverse problems in scalar transport equations, RANS turbulence models, and patient-specific cardiovascular flows, and the results have demonstrated the effectiveness and merit of the proposed method in terms of both efficiency and accuracy.

\section*{Acknowledgment}
The authors would like to acknowledge the startup funds from the College of Engineering at University of Notre Dame and funds from National Science Foundation (NSF contract CMMI-1934300) in supporting this study. We also gratefully acknowledge the discussion and collaborations from Dr. Xueyu Zhu in the early stage of this research. 

\appendix
\section{Bi-fidelity iterative ensemble Kalman method algorithm \label{Appendix:ImportantPointSelection}}
\begin{algorithm}[H]
	\footnotesize
\SetAlgoLined
\Begin{
 $V^L =[\mathbf{v}^L(\mathbf{z}_1),...,\mathbf{v}^L(\mathbf{z}_M)] $ \newline
 $\mathbf{d}[i] = (\mathbf{v}_i^L)^T\mathbf{v}_i^L$ for $i = 1,...,M$ \newline
 $P =$ \texttt{zeros}$(m,1)$, $L =$ \texttt{zeros}$(m,M)$, $k = 1$\newline 
 \While{$k \leq m$}{
  \textbf{1}. $P[k] =\arg\max \mathbf{d}[k:end]$;\newline
  \textbf{2}. Exchange $V^L[:,k]$ and $V^L[:,P[k]]$; Exchange $L[:,k]$ and $L[:,P[k]]$; Exchange $\mathbf{d}[k]$ and $\mathbf{d}[P[k]]$;\newline
  \textbf{3}. $\mathbf{r}(t) = (V^L(:,t))^TV^L(:,k) - \sum_{j=1}^{k-1}L(t,j)L(k,j)$ for $t=k+1,...,M$; \newline
  \textbf{4}. $L[k,k] = \sqrt{\mathbf{d}[k]}$;\newline
  \textbf{5}. $L[t,k] = \mathbf{r}[t]/L[k,k]$ for $t = k+1,...,M$;\newline
  \textbf{6}. $\mathbf{d}[t] = \mathbf{d}[t] - L^2[t,k]$ for $t = k+1,...,M$;\newline
  \textbf{7}. $k = k + 1$
 }
 $\gamma[:,t] = \Gamma[:,P[t]]$ for $t = 1,...,m$;\newline
 Form the truncated Gramian $G^L = LL^T$
}
\caption{Offline: important point selection for HF basis \cite{gao2019bi}}
\label{alg: pointselect}
\end{algorithm}

\begin{algorithm}[H]
	\footnotesize
	\SetAlgoLined
	{\Begin{
			Draw $\gamma_0^{(n_s)}$ from $p(\mathbf{z})$, k=1.
			\newline
			\While{$k \leq n_{iter}$}{
				\textbf{1}. $V_k^L(\gamma_k^{(n_s)})=\mathcal{F}^L(\gamma_k^{(n_s)})$;\newline
				\textbf{2}. $\mathbf{c}_k=((G^L)^{-1})^TV_k^L$;\newline
				\textbf{2}. $V^B(\gamma_k^{(n_s)})=\sum_{i=1}^m c_i \mathbf{v}^H(\mathbf{z}_i)$; \newline
				\textbf{3}. $X^B_k=[V_k^B, \gamma_k^{(n_s)}]$;
				\newline
				\textbf{3}. $\hat{X}_k^B=X_k^B+P_m^k H^T (H P_m^k H^T + P_d)^{-1}(\mathbf{y}-H X_k^B)$;
				\newline
				\textbf{4}.  $\gamma_{k+1}^{(n_s)}=\hat{\gamma}_{k}^{(n_s)}$;
				\newline
				\textbf{5}. $k = k + 1$;\newline
				\uIf{$\norm{\sqrt{P_d}(\mathbf{y}-H X_k^B)} \leq \sigma_d$}{break;}
				
			}
		}
	}
	\caption{Online: BF iterative ensemble Kalman inversion}
	\label{alg: online}
\end{algorithm}

\section{Error estimate of the BF reconstruction\label{Appendix:BFError}}
The BF construction in each case is guided by the empirical error estimation \emph{a priori}~\cite{gao2019bi}. 
\begin{figure}[H]
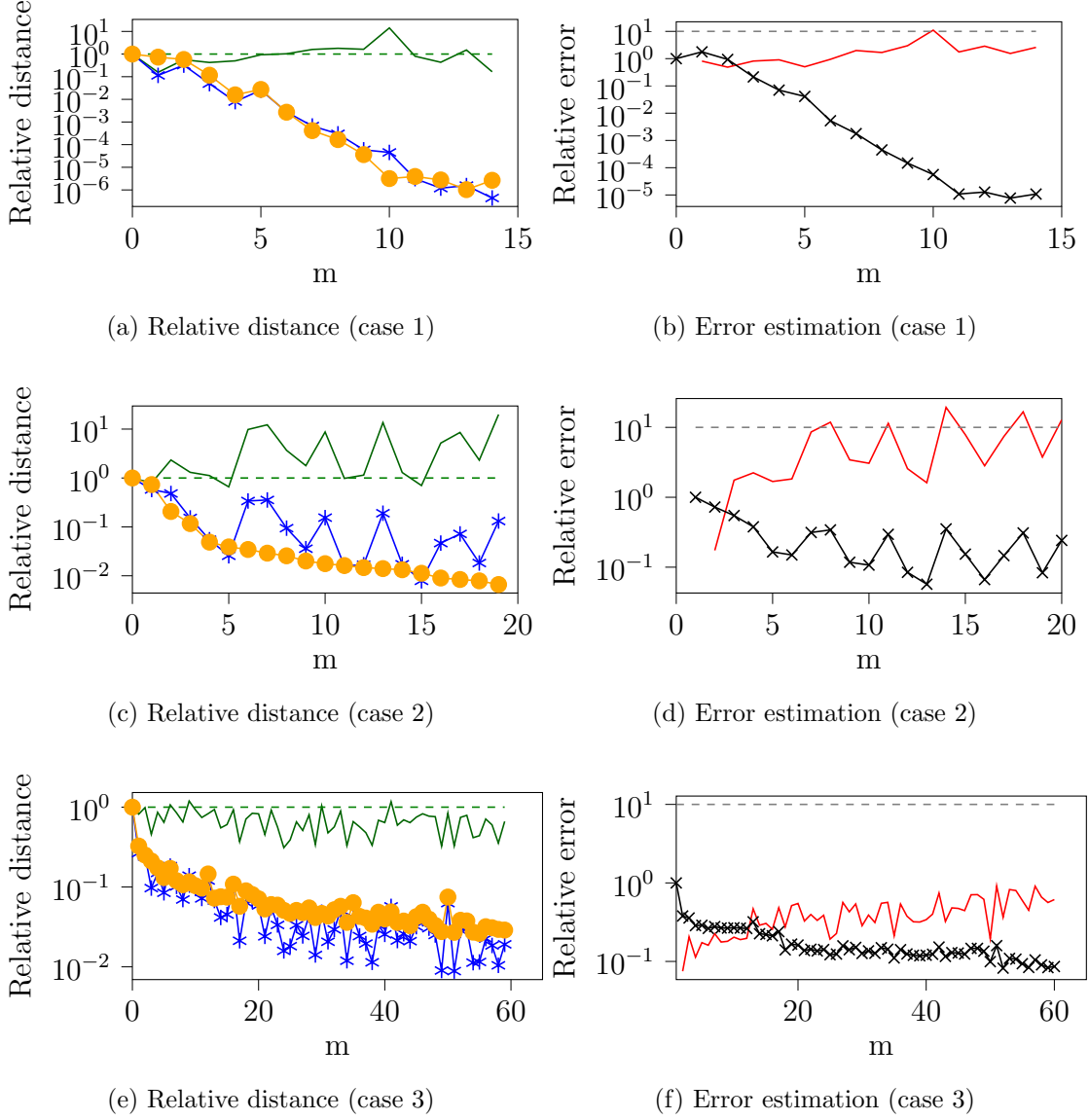

	\centering
	\subfloat[Relative distance (case 1)]
	{\includegraphics[width=0.45\textwidth,height=0.25\textwidth]{./distancePlot_case1.tikz}}
	\subfloat[Error estimation (case 1)]
	{\includegraphics[width=0.45\textwidth,height=0.235\textwidth]{./errorEstimation_case1.tikz}}\\
	
	\subfloat[Relative distance (case 2)]{\includegraphics[width=0.45\textwidth,height=0.25\textwidth]{./distancePlot_case2.tikz}}
	\subfloat[Error estimation (case 2)]{\includegraphics[width=0.45\textwidth,height=0.235\textwidth]{./errorEstimation_case2.tikz}}\\
	
	\subfloat[Relative distance (case 3)]{\includegraphics[width=0.45\textwidth,height=0.25\textwidth]{./distancePlot_case3.tikz}}
	\subfloat[Error estimation (case 3)]{\includegraphics[width=0.45\textwidth,height=0.235\textwidth]{./errorEstimation_case3.tikz}}
	\caption{The relative distance and error estimation of the BF models with respect to the number ($m$) of HF simulations used for BF reconstruction in cases 1, 2 and 3: $R_s$ (\ref{plot:Rs}), $\frac{d^H({v}^H({z}),{U}^H(\gamma^{k}))}{||{v}^H({z})||}$  (\ref{plot:distHF}), $\frac{d^L({v}^L({z}),{U}^L(\gamma^{k}))}{||{v}^L({z})||}$ (\ref{plot:distLF}), $R_e$ (\ref{plot:Re}), BF error estimation (\ref{plot:ErrorEstimation}). }
	\label{fig:case123ErrorEstimation}
\end{figure}
\noindent The model similarity metrics ($R_s$ \ref{plot:Rs}) vs. the number of HF samples are shown in the first column of Fig.~\ref{fig:case123ErrorEstimation}, where the relative distances from HF and LF solutions to the subspace spanned by the previous selected HF solutions are plotted as well. It can be seen that the model similarity metrics $R_s$ fluctuates around one, and the relative distances of both HF/LF solutions decrease at the number of important points increase, indicating that the LF models are informative for all three cases. The histories of error estimates (\ref{plot:ErrorEstimation}) and error component ratio ($R_e$ \ref{plot:Re}) are plotted in the second column of Fig.~\ref{fig:case123ErrorEstimation}. The error component ratios $R_e$ in cases 1 \& 2 almost reach to the threshold value of 10 after collecting a few HF samples, indicating the max-distance based point selection can stop. Although the $R_e$ in case 3 is far from the threshold, the estimated error does not notably decrease after 55 HF samples were collected. Therefore, the number of HF training samples in each case can be justified based on this \emph{a priori} error estimation analysis.


\begin{thebibliography}{10}
	\expandafter\ifx\csname url\endcsname\relax
	\def\url#1{\texttt{#1}}\fi
	\expandafter\ifx\csname urlprefix\endcsname\relax\def\urlprefix{URL }\fi
	\expandafter\ifx\csname href\endcsname\relax
	\def\href#1#2{#2} \def\path#1{#1}\fi
	
	\bibitem{catalano2003numerical}
	P.~Catalano, M.~Wang, G.~Iaccarino, P.~Moin, Numerical simulation of the flow
	around a circular cylinder at high reynolds numbers, International journal of
	heat and fluid flow 24~(4) (2003) 463--469.
	
	\bibitem{iaccarino2003immersed}
	G.~Iaccarino, R.~Verzicco, Immersed boundary technique for turbulent flow
	simulations, Appl. Mech. Rev. 56~(3) (2003) 331--347.
	
	\bibitem{kalitzin2005near}
	G.~Kalitzin, G.~Medic, G.~Iaccarino, P.~Durbin, Near-wall behavior of rans
	turbulence models and implications for wall functions, Journal of
	Computational Physics 204~(1) (2005) 265--291.
	
	\bibitem{evans2013isogeometric1}
	J.~A. Evans, T.~J. Hughes, Isogeometric divergence-conforming b-splines for the
	steady navier--stokes equations, Mathematical Models and Methods in Applied
	Sciences 23~(08) (2013) 1421--1478.
	
	\bibitem{evans2013isogeometric2}
	J.~A. Evans, T.~J. Hughes, Isogeometric divergence-conforming b-splines for the
	unsteady navier--stokes equations, Journal of Computational Physics 241
	(2013) 141--167.
	
	\bibitem{luo2008immersed}
	H.~Luo, R.~Mittal, X.~Zheng, S.~A. Bielamowicz, R.~J. Walsh, J.~K. Hahn, An
	immersed-boundary method for flow--structure interaction in biological
	systems with application to phonation, Journal of computational physics
	227~(22) (2008) 9303--9332.
	
	\bibitem{evans2009n}
	J.~A. Evans, Y.~Bazilevs, I.~Babu{\v{s}}ka, T.~J. Hughes, n-widths, sup--infs,
	and optimality ratios for the k-version of the isogeometric finite element
	method, Computer Methods in Applied Mechanics and Engineering 198~(21-26)
	(2009) 1726--1741.
	
	\bibitem{evans2013isogeometric3}
	J.~A. Evans, T.~J. Hughes, Isogeometric divergence-conforming b-splines for the
	darcy--stokes--brinkman equations, Mathematical Models and Methods in Applied
	Sciences 23~(04) (2013) 671--741.
	
	\bibitem{zheng2009computational}
	X.~Zheng, S.~Bielamowicz, H.~Luo, R.~Mittal, A computational study of the
	effect of false vocal folds on glottal flow and vocal fold vibration during
	phonation, Annals of biomedical engineering 37~(3) (2009) 625--642.
	
	\bibitem{plessix2006review}
	R.-E. Plessix, A review of the adjoint-state method for computing the gradient
	of a functional with geophysical applications, Geophysical Journal
	International 167~(2) (2006) 495--503.
	
	\bibitem{giles2000introduction}
	M.~B. Giles, N.~A. Pierce, An introduction to the adjoint approach to design,
	Flow, turbulence and combustion 65~(3-4) (2000) 393--415.
	
	\bibitem{dow2011quantification}
	E.~Dow, Q.~Wang, Quantification of structural uncertainties in the $k-\omega$
	turbulence model, in: 52nd AIAA/ASME/ASCE/AHS/ASC Structures, Structural
	Dynamics and Materials Conference 19th AIAA/ASME/AHS Adaptive Structures
	Conference 13t, 2011, p. 1762.
	
	\bibitem{talnikar2017unsteady}
	C.~Talnikar, Q.~Wang, G.~M. Laskowski, Unsteady adjoint of pressure loss for a
	fundamental transonic turbine vane, Journal of Turbomachinery 139~(3).
	
	\bibitem{pires2001tsunami}
	C.~Pires, P.~M. Miranda, Tsunami waveform inversion by adjoint methods, Journal
	of Geophysical Research: Oceans 106~(C9) (2001) 19773--19796.
	
	\bibitem{wang2019discrete}
	M.~Wang, Q.~Wang, T.~A. Zaki, Discrete adjoint of fractional-step
	incompressible navier-stokes solver in curvilinear coordinates and
	application to data assimilation, Journal of Computational Physics 396 (2019)
	427--450.
	
	\bibitem{rangarajan2020adjoint}
	A.~Rangarajan, G.~May, V.~Dolejsi, Adjoint-based anisotropic $hp$-adaptation
	for discontinuous {Galerkin} methods using a continuous mesh model, Journal
	of Computational Physics (2020) 109321.
	
	\bibitem{wang2020identification}
	Z.~Wang, X.~Huan, K.~Garikipati, Identification of the partial differential
	equations governing microstructure evolution in materials: Inference over
	incomplete, sparse and spatially non-overlapping data, arXiv preprint
	arXiv:2001.04816.
	
	\bibitem{nielsen2010discrete}
	E.~J. Nielsen, B.~Diskin, N.~K. Yamaleev, Discrete adjoint-based design
	optimization of unsteady turbulent flows on dynamic unstructured grids, AIAA
	journal 48~(6) (2010) 1195--1206.
	
	\bibitem{asch2016data}
	M.~Asch, M.~Bocquet, M.~Nodet, Data assimilation: methods, algorithms, and
	applications, Vol.~11, SIAM, 2016.
	
	\bibitem{cotter2009bayesian}
	S.~L. Cotter, M.~Dashti, J.~C. Robinson, A.~M. Stuart, Bayesian inverse
	problems for functions and applications to fluid mechanics, Inverse problems
	25~(11) (2009) 115008.
	
	\bibitem{ou2019new}
	N.~Ou, L.~Jiang, G.~Lin, A new bi-fidelity model reduction method for bayesian
	inverse problems, International Journal for Numerical Methods in Engineering
	119~(10) (2019) 941--963.
	
	\bibitem{zhang2019efficient}
	J.~Zhang, M.~D. Shields, Efficient monte carlo resampling for probability
	measure changes from bayesian updating, Probabilistic Engineering Mechanics
	55 (2019) 54--66.
	
	\bibitem{zhang2018quantification}
	J.~Zhang, M.~D. Shields, On the quantification and efficient propagation of
	imprecise probabilities resulting from small datasets, Mechanical Systems and
	Signal Processing 98 (2018) 465--483.
	
	\bibitem{zhang2018effect}
	J.~Zhang, M.~D. Shields, The effect of prior probabilities on quantification
	and propagation of imprecise probabilities resulting from small datasets,
	Computer Methods in Applied Mechanics and Engineering 334 (2018) 483--506.
	
	\bibitem{li2015adaptive}
	W.~Li, G.~Lin, An adaptive importance sampling algorithm for bayesian inversion
	with multimodal distributions, Journal of Computational Physics 294 (2015)
	173--190.
	
	\bibitem{morzfeld2018iterative}
	M.~Morzfeld, M.~S. Day, R.~W. Grout, G.~S. Heng~Pau, S.~A. Finsterle, J.~B.
	Bell, Iterative importance sampling algorithms for parameter estimation, SIAM
	Journal on Scientific Computing 40~(2) (2018) B329--B352.
	
	\bibitem{uzun2019structural}
	M.~Uzun, H.~Sun, D.~Smit, O.~B{\"u}y{\"u}k{\"o}zt{\"u}rk, Structural damage
	detection using bayesian inference and seismic interferometry, Structural
	Control and Health Monitoring 26~(11) (2019) e2445.
	
	\bibitem{iglesias2013ensemble}
	M.~A. Iglesias, K.~J. Law, A.~M. Stuart, Ensemble kalman methods for inverse
	problems, Inverse Problems 29~(4) (2013) 045001.
	
	\bibitem{liu2016nonlocal}
	Y.~Liu, W.~Sun, Z.~Yuan, J.~Fish, A nonlocal multiscale discrete-continuum
	model for predicting mechanical behavior of granular materials, International
	Journal for Numerical Methods in Engineering 106~(2) (2016) 129--160.
	
	\bibitem{sun2015statistical}
	H.~Sun, D.~Feng, Y.~Liu, M.~Q. Feng, Statistical regularization for
	identification of structural parameters and external loadings using state
	space models, Computer-Aided Civil and Infrastructure Engineering 30~(11)
	(2015) 843--858.
	
	\bibitem{liu2016determining}
	Y.~Liu, W.~Sun, J.~Fish, Determining material parameters for critical state
	plasticity models based on multilevel extended digital database, Journal of
	Applied Mechanics 83~(1).
	
	\bibitem{sun2015hybrid}
	H.~Sun, R.~Betti, A hybrid optimization algorithm with bayesian inference for
	probabilistic model updating, Computer-Aided Civil and Infrastructure
	Engineering 30~(8) (2015) 602--619.
	
	\bibitem{sun2013identification}
	H.~Sun, H.~Lu{\c{s}}, R.~Betti, Identification of structural models using a
	modified artificial bee colony algorithm, Computers \& Structures 116 (2013)
	59--74.
	
	\bibitem{feng2015simultaneous}
	D.~Feng, H.~Sun, M.~Q. Feng, Simultaneous identification of bridge structural
	parameters and vehicle loads, Computers \& Structures 157 (2015) 76--88.
	
	\bibitem{sun2016probabilistic}
	H.~Sun, O.~B{\"u}y{\"u}k{\"o}zt{\"u}rk, Probabilistic updating of building
	models using incomplete modal data, Mechanical Systems and Signal Processing
	75 (2016) 27--40.
	
	\bibitem{evensen2003ensemble}
	G.~Evensen, The ensemble {Kalman} filter: Theoretical formulation and practical
	implementation, Ocean dynamics 53~(4) (2003) 343--367.
	
	\bibitem{evensen2009data}
	G.~Evensen, Data assimilation: the ensemble Kalman filter, Springer Science \&
	Business Media, 2009.
	
	\bibitem{carrassi2018data}
	A.~Carrassi, M.~Bocquet, L.~Bertino, G.~Evensen, Data assimilation in the
	geosciences: An overview of methods, issues, and perspectives, Wiley
	Interdisciplinary Reviews: Climate Change 9~(5) (2018) e535.
	
	\bibitem{schillings2017analysis}
	C.~Schillings, A.~M. Stuart, Analysis of the ensemble {Kalman} filter for
	inverse problems, SIAM Journal on Numerical Analysis 55~(3) (2017)
	1264--1290.
	
	\bibitem{chen2012ensemble}
	Y.~Chen, D.~S. Oliver, Ensemble randomized maximum likelihood method as an
	iterative ensemble smoother, Mathematical Geosciences 44~(1) (2012) 1--26.
	
	\bibitem{chen2013levenberg}
	Y.~Chen, D.~S. Oliver, Levenberg--marquardt forms of the iterative ensemble
	smoother for efficient history matching and uncertainty quantification,
	Computational Geosciences 17~(4) (2013) 689--703.
	
	\bibitem{schillings2018convergence}
	C.~Schillings, A.~M. Stuart, Convergence analysis of ensemble kalman inversion:
	the linear, noisy case, Applicable Analysis 97~(1) (2018) 107--123.
	
	\bibitem{evensen2018analysis}
	G.~Evensen, Analysis of iterative ensemble smoothers for solving inverse
	problems, Computational Geosciences 22~(3) (2018) 885--908.
	
	\bibitem{blomker2019well}
	D.~Bl{\"o}mker, C.~Schillings, P.~Wacker, S.~Weissmann, Well posedness and
	convergence analysis of the ensemble {Kalman} inversion, Inverse Problems
	35~(8) (2019) 085007.
	
	\bibitem{wu2019improving}
	J.~Wu, J.-X. Wang, S.~C. Shadden, Improving the convergence of the iterative
	ensemble kalman filter by resampling, arXiv preprint arXiv:1910.04247.
	
	\bibitem{iglesias2016regularizing}
	M.~A. Iglesias, A regularizing iterative ensemble kalman method for
	pde-constrained inverse problems, Inverse Problems 32~(2) (2016) 025002.
	
	\bibitem{wu2019adding}
	J.~Wu, J.-X. Wang, S.~C. Shadden, Adding constraints to bayesian inverse
	problems, in: Proceedings of the AAAI Conference on Artificial Intelligence,
	Vol.~33, 2019, pp. 1666--1673.
	
	\bibitem{zhang2019regularization}
	X.-L. Zhang, C.~Michel{\'e}n-Str{\"o}fer, H.~Xiao, Regularization of ensemble
	kalman methods for inverse problems, arXiv preprint arXiv:1910.01292.
	
	\bibitem{albers2019ensemble}
	D.~J. Albers, P.-A. Blancquart, M.~E. Levine, E.~E. Seylabi, A.~Stuart,
	Ensemble kalman methods with constraints, Inverse Problems 35~(9) (2019)
	095007.
	
	\bibitem{chada2019incorporation}
	N.~K. Chada, C.~Schillings, S.~Weissmann, On the incorporation of
	box-constraints for ensemble kalman inversion, arXiv preprint
	arXiv:1908.00696.
	
	\bibitem{iglesias2015iterative}
	M.~A. Iglesias, Iterative regularization for ensemble data assimilation in
	reservoir models, Computational Geosciences 19~(1) (2015) 177--212.
	
	\bibitem{wang2016data}
	J.-X. Wang, H.~Xiao, Data-driven cfd modeling of turbulent flows through
	complex structures, International Journal of Heat and Fluid Flow 62 (2016)
	138--149.
	
	\bibitem{wang2018inferring}
	J.-X. Wang, H.~Tang, H.~Xiao, R.~Weiss, Inferring tsunami flow depth and flow
	speed from sediment deposits based on ensemble kalman filtering, Geophysical
	Journal International 212~(1) (2018) 646--658.
	
	\bibitem{tang2018tsuflind}
	H.~Tang, J.~Wang, R.~Weiss, H.~Xiao, Tsuflind-enkf: Inversion of tsunami flow
	depth and flow speed from deposits with quantified uncertainties, Marine
	Geology 396 (2018) 16--25.
	
	\bibitem{kato2013approach}
	H.~Kato, S.~Obayashi, Approach for uncertainty of turbulence modeling based on
	data assimilation technique, Computers \& Fluids 85 (2013) 2--7.
	
	\bibitem{xiao2016quantifying}
	H.~Xiao, J.-L. Wu, J.-X. Wang, R.~Sun, C.~Roy, Quantifying and reducing
	model-form uncertainties in reynolds-averaged navier--stokes simulations: A
	data-driven, physics-informed bayesian approach, Journal of Computational
	Physics 324 (2016) 115--136.
	
	\bibitem{wu2016bayesian}
	J.-L. Wu, J.-X. Wang, H.~Xiao, A bayesian calibration--prediction method for
	reducing model-form uncertainties with application in rans simulations, Flow,
	Turbulence and Combustion 97~(3) (2016) 761--786.
	
	\bibitem{mons2016reconstruction}
	V.~Mons, J.-C. Chassaing, T.~Gomez, P.~Sagaut, Reconstruction of unsteady
	viscous flows using data assimilation schemes, Journal of Computational
	Physics 316 (2016) 255--280.
	
	\bibitem{arnold2017uncertainty}
	A.~Arnold, C.~Battista, D.~Bia, Y.~Z. German, R.~L. Armentano, H.~Tran, M.~S.
	Olufsen, Uncertainty quantification in a patient-specific one-dimensional
	arterial network model: Enkf-based inflow estimator, Journal of Verification,
	Validation and Uncertainty Quantification 2~(1).
	
	\bibitem{lal2017non}
	R.~Lal, F.~Nicoud, E.~Le~Bars, J.~Deverdun, F.~Molino, V.~Costalat,
	B.~Mohammadi, Non invasive blood flow features estimation in cerebral
	arteries from uncertain medical data, Annals of biomedical engineering
	45~(11) (2017) 2574--2591.
	
	\bibitem{wang2019data}
	J.-X. Wang, X.~Hu, S.~C. Shadden, Data-augmented modeling of intracranial
	pressure, Annals of biomedical engineering 47~(3) (2019) 714--730.
	
	\bibitem{zoccarato2016data}
	C.~Zoccarato, D.~Ba{\`u}, M.~Ferronato, G.~Gambolati, A.~Alzraiee, P.~Teatini,
	Data assimilation of surface displacements to improve geomechanical
	parameters of gas storage reservoirs, Journal of Geophysical Research: Solid
	Earth 121~(3) (2016) 1441--1461.
	
	\bibitem{iglesias2018ensemble}
	M.~Iglesias, Z.~Sawlan, M.~Scavino, R.~Tempone, C.~Wood, Ensemble-marginalized
	kalman filter for linear time-dependent pdes with noisy boundary conditions:
	Application to heat transfer in building walls, Inverse Problems 34~(7)
	(2018) 075008.
	
	\bibitem{sousa2019computational}
	J.~Sousa, C.~Gorl{\'e}, Computational urban flow predictions with bayesian
	inference: Validation with field data, Building and Environment 154 (2019)
	13--22.
	
	\bibitem{yang2020non}
	F.-l. Yang, L.~Yan, A non-intrusive reduced basis eki for time fractional
	diffusion inverse problems, Acta Mathematicae Applicatae Sinica, English
	Series 36~(1) (2020) 183--202.
	
	\bibitem{narayan2014stochastic}
	A.~Narayan, C.~Gittelson, D.~Xiu, A stochastic collocation algorithm with
	multifidelity models, SIAM Journal on Scientific Computing 36~(2) (2014)
	A495--A521.
	
	\bibitem{zhu2014computational}
	X.~Zhu, A.~Narayan, D.~Xiu, Computational aspects of stochastic collocation
	with multifidelity models, SIAM/ASA Journal on Uncertainty Quantification
	2~(1) (2014) 444--463.
	
	\bibitem{hoel2016multilevel}
	H.~Hoel, K.~J. Law, R.~Tempone, Multilevel ensemble kalman filtering, SIAM
	Journal on Numerical Analysis 54~(3) (2016) 1813--1839.
	
	\bibitem{narayan2012stochastic}
	A.~Narayan, D.~Xiu, Stochastic collocation methods on unstructured grids in
	high dimensions via interpolation, SIAM Journal on Scientific Computing
	34~(3) (2012) A1729--A1752.
	
	\bibitem{gao2019bi}
	H.~Gao, X.~Zhu, J.-X. Wang, A bi-fidelity surrogate modeling approach for
	uncertainty propagation in three-dimensional hemodynamic simulations, arXiv
	preprint arXiv:1908.10197.
	
	\bibitem{hampton2018practical}
	J.~Hampton, H.~R. Fairbanks, A.~Narayan, A.~Doostan, Practical error bounds for
	a non-intrusive bi-fidelity approach to parametric/stochastic model
	reduction, Journal of Computational Physics 368 (2018) 315--332.
	
	\bibitem{skinner2019reduced}
	R.~W. Skinner, A.~Doostan, E.~L. Peters, J.~A. Evans, K.~E. Jansen,
	Reduced-basis multifidelity approach for efficient parametric study of naca
	airfoils, AIAA Journal 57~(4) (2019) 1481--1491.
	
	\bibitem{fairbanks2020bi}
	H.~R. Fairbanks, L.~Jofre, G.~Geraci, G.~Iaccarino, A.~Doostan, Bi-fidelity
	approximation for uncertainty quantification and sensitivity analysis of
	irradiated particle-laden turbulence, Journal of Computational Physics 402
	(2020) 108996.
	
	\bibitem{stein1987large}
	M.~Stein, Large sample properties of simulations using latin hypercube
	sampling, Technometrics 29~(2) (1987) 143--151.
	
	\bibitem{shields2016generalization}
	M.~D. Shields, J.~Zhang, The generalization of latin hypercube sampling,
	Reliability Engineering \& System Safety 148 (2016) 96--108.
	
	\bibitem{caretto1973two}
	L.~Caretto, A.~Gosman, S.~Patankar, D.~Spalding, Two calculation procedures for
	steady, three-dimensional flows with recirculation, in: Proceedings of the
	third international conference on numerical methods in fluid mechanics,
	Springer, 1973, pp. 60--68.
	
	\bibitem{rhie1983numerical}
	C.~Rhie, W.~L. Chow, Numerical study of the turbulent flow past an airfoil with
	trailing edge separation, AIAA journal 21~(11) (1983) 1525--1532.
	
	\bibitem{canuto2019uncertainty}
	C.~Canuto, S.~Pieraccini, D.~Xiu, Uncertainty quantification of discontinuous
	outputs via a non-intrusive bifidelity strategy, Journal of Computational
	Physics 398 (2019) 108885.
	
	\bibitem{durbin2018some}
	P.~A. Durbin, Some recent developments in turbulence closure modeling, Annual
	Review of Fluid Mechanics 50 (2018) 77--103.
	
	\bibitem{spalart2009detached}
	P.~R. Spalart, Detached-eddy simulation, Annual review of fluid mechanics 41
	(2009) 181--202.
	
	\bibitem{wang2016incorporating}
	J.~Wang, J.-L. Wu, H.~Xiao, Incorporating prior knowledge for quantifying and
	reducing model-form uncertainty in rans simulations, International Journal
	for Uncertainty Quantification 6~(2).
	
	\bibitem{edeling2018bayesian}
	W.~N. Edeling, M.~Schmelzer, R.~P. Dwight, P.~Cinnella, Bayesian predictions of
	reynolds-averaged navier--stokes uncertainties using maximum a posteriori
	estimates, AIAA Journal 56~(5) (2018) 2018--2029.
	
	\bibitem{wang2017physics}
	J.-X. Wang, J.-L. Wu, H.~Xiao, Physics-informed machine learning approach for
	reconstructing reynolds stress modeling discrepancies based on dns data,
	Physical Review Fluids 2~(3) (2017) 034603.
	
	\bibitem{ling2016machine}
	J.~Ling, R.~Jones, J.~Templeton, Machine learning strategies for systems with
	invariance properties, Journal of Computational Physics 318 (2016) 22--35.
	
	\bibitem{duraisamy2019turbulence}
	K.~Duraisamy, G.~Iaccarino, H.~Xiao, Turbulence modeling in the age of data,
	Annual Review of Fluid Mechanics 51 (2019) 357--377.
	
	\bibitem{wang2019prediction}
	J.-X. Wang, J.~Huang, L.~Duan, H.~Xiao, Prediction of reynolds stresses in
	high-mach-number turbulent boundary layers using physics-informed machine
	learning, Theoretical and Computational Fluid Dynamics 33~(1) (2019) 1--19.
	
	\bibitem{yang2019predictive}
	X.~Yang, S.~Zafar, J.-X. Wang, H.~Xiao, Predictive large-eddy-simulation wall
	modeling via physics-informed neural networks, Physical Review Fluids 4~(3)
	(2019) 034602.
	
	\bibitem{launder1983numerical}
	B.~E. Launder, D.~B. Spalding, The numerical computation of turbulent flows,
	in: Numerical prediction of flow, heat transfer, turbulence and combustion,
	Elsevier, 1983, pp. 96--116.
	
	\bibitem{fureby2010homogenization}
	C.~Fureby, Homogenization based les for turbulent combustion, Flow, turbulence
	and combustion 84~(3) (2010) 459--480.
	
	\bibitem{valen2018real}
	K.~Valen-Sendstad, A.~W. Bergersen, Y.~Shimogonya, L.~Goubergrits, J.~Bruening,
	J.~Pallares, S.~Cito, S.~Piskin, K.~Pekkan, A.~J. Geers, et~al., Real-world
	variability in the prediction of intracranial aneurysm wall shear stress: The
	2015 international aneurysm cfd challenge, Cardiovascular engineering and
	technology 9~(4) (2018) 544--564.
	
	\bibitem{tipping1999probabilistic}
	M.~E. Tipping, C.~M. Bishop, Probabilistic principal component analysis,
	Journal of the Royal Statistical Society: Series B (Statistical Methodology)
	61~(3) (1999) 611--622.
	
	\bibitem{jou2008wall}
	L.-D. Jou, D.~H. Lee, H.~Morsi, M.~E. Mawad, Wall shear stress on ruptured and
	unruptured intracranial aneurysms at the internal carotid artery, American
	Journal of Neuroradiology 29~(9) (2008) 1761--1767.
	
\end{thebibliography}

\end{document}